\documentclass[journal]{IEEEtran}

 \usepackage{graphicx} 
  \DeclareGraphicsExtensions{.eps}
  
  \usepackage{tikz}
  \usepackage{pgfplots,pgfplotstable}
  \usetikzlibrary{pgfplots.groupplots}
  \pgfplotsset{compat=1.5}

\usepackage{amsmath}

\usepackage{empheq} 

\usepackage{textcomp}
\usetikzlibrary{shapes,arrows}

\usepackage{stfloats} 

\usepackage[hyphens]{url}
\usepackage[hidelinks]{hyperref}
\hypersetup{colorlinks=true,linkcolor=black,citecolor=black,filecolor=black,urlcolor=black,breaklinks=true}
\usepackage{cite} 

\usepackage{color}

\usepackage{array} 

\usepackage[export]{adjustbox} 

\usepackage[normalem]{ulem} 

\ifCLASSINFOpdf
\else
\fi
\makeatletter
\def\bstctlcite{\@ifnextchar[{\@bstctlcite}{\@bstctlcite[@auxout]}}
\def\@bstctlcite[#1]#2{\@bsphack
 \@for\@citeb:=#2\do{%
   \edef\@citeb{\expandafter\@firstofone\@citeb}%
   \if@filesw\immediate\write\csname #1\endcsname{\string\citation{\@citeb}}\fi}%
 \@esphack}
\makeatother

\begin{document}
%

\bstctlcite{IEEEexample:BSTcontrol}

\title{Simultaneous Scheduling of Multiple Frequency Services in Stochastic Unit Commitment}
%
%
%

\author{Luis~Badesa,~\IEEEmembership{Student Member,~IEEE,}
        Fei~Teng,~\IEEEmembership{Member,~IEEE,}
        and~Goran~Strbac,~\IEEEmembership{Member,~IEEE}
\thanks{A portion of this work has been supported by the Engineering and Physical Sciences Research Council under grants EP/K002252/1 and EP/R045518/1.}
\thanks{
The authors are with the Department of Electrical and Electronic Engineering, Imperial College London, SW7 2AZ London, U.K. (email: luis.badesa@imperial.ac.uk).}
}

%
%

\markboth{IEEE Transactions on Power Systems, March~2019}%
{Shell \MakeLowercase{\textit{et al.}}: Bare Demo of IEEEtran.cls for IEEE Journals}
%



\maketitle

\begin{abstract}
The reduced level of system inertia in low-carbon power grids increases the need for alternative frequency services. However, simultaneously optimising the provision of these services in the scheduling process, subject to significant uncertainty, is a complex task given the challenge of linking the steady-state optimisation with frequency dynamics. This paper proposes a novel frequency-constrained Stochastic Unit Commitment (SUC) model which, for the first time, co-optimises energy production along with the provision of synchronised and synthetic inertia, Enhanced Frequency Response (EFR), Primary Frequency Response (PFR) and a dynamically-reduced largest power infeed. The contribution of load damping is modelled through a linear inner approximation. The effectiveness of the proposed model is demonstrated through several case studies for Great Britain's 2030 power system, which highlight the synergies and conflicts among alternative frequency services, as well as the significant economic savings and carbon reduction achieved by simultaneously optimising all these services. 
\end{abstract}


%
\IEEEpeerreviewmaketitle

\section*{Nomenclature}
\addcontentsline{toc}{section}{Nomenclature}

\subsection*{Indices and Sets}
\begin{IEEEdescription}[\IEEEusemathlabelsep\IEEEsetlabelwidth{$\textrm{RoCoF}_{\textrm{ma}}$}]
\item[$g,\,\, \mathcal{G}$] Index, Set of generators.
\item[$l,\,\, \mathcal{L}$] Index, Set of bits in binary expansion of PFR.
\item[$n,\,\, \mathcal{N}$] Index, Set of nodes in the scenario tree.
\item[$p,\,\, \mathcal{P}$] Index, Set of overestimating planes.
\item[$s,\,\, \mathcal{S}$] Index, Set of storage units.
\end{IEEEdescription}

\subsection*{Constants and Variables \normalfont{(not related to the optimisation)}}
\begin{IEEEdescription}[\IEEEusemathlabelsep\IEEEsetlabelwidth{$\textrm{RoCoF}_{\textrm{ma}}$}]
\item[$\tau_\textrm{b}$] Time-constant of the BESS dynamics (s).
\item[$\tau_\textrm{g}$] Time-constant of the generator dynamics (s).
\item[$\textrm{K}_\textrm{b}$] Droop gain for the BESS dynamics (Hz/MW).
\item[$\textrm{K}_\textrm{g}$] Droop gain for the generator dynamics (Hz/MW).
\item[$t^*$] Time when the frequency nadir is reached (s).
\end{IEEEdescription}

\subsection*{Constants}
\begin{IEEEdescription}[\IEEEusemathlabelsep\IEEEsetlabelwidth{$\textrm{RoCoF}_{\textrm{ma}}$}]
\item[$\Delta\tau(n)$] Time-step corresponding to node $n$ (h).
\item[$\Delta f_{\textrm{max}}$] Maximum admissible frequency deviation (Hz).
\vspace{-10pt}
\item[$\Delta f^{\textrm{ss}}_{\textrm{max}}$] Maximum frequency deviation at quasi-steady-state (Hz).
\item[$\pi(n)$] Probability of reaching node $n$.
\vspace{1pt}
\item[$\textrm{a}_p,\,\textrm{b}_p,\,\textrm{c}_p$] Parameters of overestimating plane $p$.
\vspace{2pt}
\item[$\textrm{c}^{\textrm{m}}_g$] Marginal cost of generating units $g$ (\pounds/MWh).
\vspace{3pt}
\item[$\textrm{c}^{\textrm{nl}}_g$] No-load cost of generating units $g$ (\pounds/h).
\vspace{3pt}
\item[$\textrm{c}^{\textrm{st}}_g$] Startup cost of generating units $g$ (\pounds).
\vspace{2pt}
\item[$\textrm{D}$] Load-damping factor (\%/Hz).
\item[$\textrm{D}'$] Multiplication of $\textrm{D} \cdot \textrm{P}_{\textrm{D}}$ (MW/Hz).
\item[$f_0$] Nominal frequency of the power grid (Hz).
\item[$\textrm{H}_g$] Inertia constant of generating units $g$ (s).
\item[$\textrm{H}_\textrm{L}$] Inertia constant of generator producing $P_{\textrm{L}}$ (s).
\vspace{1pt}
\item[$\textrm{H}_\textrm{W}$] Inertia constant of wind turbines (s).
\item[$\textrm{N}_\mathcal{M}$] Number of must-run generators.
\item[$\textrm{P}_{\textrm{D}}$] Total demand (MW).
\vspace{1pt}
\item[$\textrm{P}_{\textrm{g}}^{\textrm{max}}$] Maximum power output of units $g$ (MW).
\vspace{3pt}
\item[$\textrm{P}_{\textrm{L}}^{\textrm{max}}$] Upper bound for the largest power infeed (MW).
\item[$\textrm{RoCoF}_{\textrm{max}}$] Maximum admissible RoCoF (Hz/s).
\item[$R_{\mathcal{S},\textrm{max}}$] Upper bound for $R_{\mathcal{S}}$ (MW).
\item[$\textrm{T}_{\textrm{g}}$] Delivery time of PFR (s).
\item[$\textrm{T}_{\textrm{s}}$] Delivery time of EFR (s).
\end{IEEEdescription}


\subsection*{Decision Variables \normalfont{(continuous unless otherwise indicated)}}
\begin{IEEEdescription}[\IEEEusemathlabelsep\IEEEsetlabelwidth{$\textrm{RoCoF}_{\textrm{ma}}$}]
\item[$k_l$] Aux. variables for linearising $R_\mathcal{S}\cdot z_l$ (MW).
\item[$m_l$] Aux. variables for linearising $H\cdot z_l$ (MW$\cdot \textrm{s}$).
\vspace{2pt}
\item[$N_g^\textrm{sg}(n)$] Number of units $g$ that start generating in node $n$.
\item[$N_g^{\textrm{up}}$] Number of online generating units of type $g$.
\vspace{1pt}
\item[$P_g$] Power produced by generating units $g$ (MW).
\vspace{1pt}
\item[$P_{\textrm{L}}$] Largest power infeed (MW).
\item[$P_\mathcal{M}$] Power produced by must-run units (MW).
\item[$P_\textrm{W}$] Online power from wind (MW).
\item[$R_g$] PFR provision from generating units $g$ (MW).
\item[$R_s$] EFR provision from storage units $s$ (MW).
\item[$z_l$] Binary variables for binary expansion of PFR.
\end{IEEEdescription}

\subsection*{Linear Expressions \normalfont{(linear combinations of decision variables)}}
\begin{IEEEdescription}[\IEEEusemathlabelsep\IEEEsetlabelwidth{$\textrm{RoCoF}_{\textrm{ma}}$}]
\item[$C_g(n)$] Operating cost of units $g$ at node $n$ (\pounds).
\item[$H$] System inertia after the loss of $P_{\textrm{L}}$ (MW$\cdot \textrm{s}$).
\item[$R_{\mathcal{G}}$] Total PFR from all generators (MW).
\item[$R_{\mathcal{S}}$] Total EFR from all storage units (MW).
\end{IEEEdescription}

\section{Introduction}
%
%
%
%
\IEEEPARstart{W}{ith} the increasing penetration of nonsynchronous renewable energy sources (RES), post-fault frequency security is becoming a greater concern. Maintaining frequency security consists on keeping the system frequency within safe boundaries. In the event of a generation loss, the subsequent frequency excursion is contained by frequency services as inertia, load damping and various forms of frequency response (FR). As nonsynchronous RES displace thermal units but currently do not provide inertia, the level of system inertia in low-carbon grids is greatly reduced, leading to a higher risk of frequency instability. 


 

Traditionally, the scheduling algorithms for power grids only enforced the total amount of FR to be above a pre-defined threshold \cite{UCRestrepoGaliana}. Due to the abundant system inertia from conventional thermal plants, this threshold was mainly driven by the requirement for frequency to return to its nominal value following a power outage, i.e. the frequency steady-state limit. However, given that system inertia may be highly scarce at times when high RES production coincides with low demand, the limit for frequency nadir during the transient period starts to drive the FR requirement, which requires explicit consideration of post-fault frequency dynamics in order to maintain system security. 

In this context, the research community is increasingly focused on implementing various frequency-security constraints in Optimal Power Flow and Unit Commitment (UC). The key challenge lies on the mathematical complexity of incorporating the differential-equation-driven frequency evolution into the algebraic-equation-constrained optimisation problem. Diverse ways of obtaining the frequency-security constraints, namely maximum admissible Rate of Change of Frequency (RoCoF), minimum admissible frequency at nadir and minimum admissible frequency at quasi-steady-state, have been discussed in \cite{UCRestrepoGaliana,OPFChavez,ElaI,LinearizedUC,FeiStochastic,DartmouthUC,UCFaroe,UC_EFRpaper,IowaThesis,ERCOT_EFR}. Certain assumptions were made in these studies to overcome the mathematical difficulties and obtain relatively simple algebraic constraints suitable for implementation in a Mixed-Integer Linear Program (MILP). Reference \cite{UCRestrepoGaliana} focused only on the quasi-steady-state condition, while \cite{OPFChavez} considered the system inertia to be a fixed, known value. Linearisation techniques are proposed in \cite{LinearizedUC,FeiStochastic,DartmouthUC,UC_EFRpaper} to transform their deduced nonlinear constraints. Except \cite{LinearizedUC,FeiStochastic,DartmouthUC}, all previous work neglected the effect of load damping on frequency evolution, leading to an overestimation of the frequency challenge and preventing incentivising the potential damping providers.

Furthermore, the reduced level of inertia increases the need for alternative frequency services in order to achieve a cost-effective integration of RES. In the past, Primary Frequency Response (PFR) was the only service considered to contain frequency decline, but system operators start to investigate the possibility to incorporate frequency services with shorter delivery time. National Grid in Great Britain (GB) has recently introduced a new FR service, Enhanced Frequency Response (EFR) \cite{NationalGridPLossInertia}, for which response must be delivered within one second, as opposed to ten seconds for PFR. A pre-defined amount of 200MW EFR was procured throughout the year in 2017, while the optimal portfolio of frequency services actually varies along with changes in system conditions. Reference \cite{UC_EFRpaper} recently proposed a UC framework that optimises EFR provision from Battery Energy Storage Systems (BESS), relying on the assumption that EFR is delivered instantaneously after a power outage. This assumption greatly reduces the mathematical complexity of the resulting constraints, but fails to reflect the actual characteristics of EFR.

According to the security standard, FR is scheduled to cover the loss of the largest power infeed. Recently, National Grid has proposed to reduce this largest infeed \cite{NationalGridPLossInertia}, as a measure to tackle the frequency stability challenge under certain system conditions. To the best of our knowledge, the size of the largest contingency has not yet been modelled as a decision variable in any frequency-secured UC formulation. Reference \cite{OMalleyDeload} considered a variable largest power outage in a competitive market-dispatch framework, deducing the frequency-security constraints from dynamic simulation results. However, this approach only allows to accurately model a limited number of the operational conditions. 

Moreover, synthetic inertia (SI) provision from wind turbines has been considered as an alternative to resolve the frequency-decline challenge. Although some studies have analysed the potential value of SI such as \cite{SyntheticEla}, it is not clear yet how SI would affect the value of alternative frequency services.
 
Given this background, this paper proposes analytical constraints for a secure post-fault frequency evolution in a Stochastic Unit Commitment (SUC) model. The key contributions of this work are three-fold:
\begin{enumerate}
	\item This paper proposes novel frequency stability constraints that, for the first time, allow to simultaneously co-optimise the provision of synchronised and synthetic inertia, PFR, EFR and a dynamically-reduced largest power infeed. The contribution of load damping to support the frequency nadir is also modelled through a linear inner approximation.
    \item The resulting nonconvex frequency-nadir constraint is efficiently linearised for implementation in a computationally-demanding SUC framework. The accuracy and computational efficiency of the linearisation method are explicitly quantified, and frequency security is guaranteed in all cases. 
    \item The proposed model is applied to the GB 2030 system to demonstrate the benefits of simultaneously co-optimising a portfolio of diverse frequency services, as well as the impact of competition among these services.
\end{enumerate}

The rest of this paper is organised as follows: Section \ref{SUCsection} gives a brief description of the SUC model. The mathematical deduction of the frequency-security constraints is presented in Section \ref{SectionFrequency}. Section \ref{SectionCaseStudies} discusses the results of several case studies, while Section \ref{SectionConclusion} concludes the paper.

\section{Stochastic Unit Commitment} \label{SUCsection}

\begin{figure}[!t]
\centering
\includegraphics[width=3.1in]{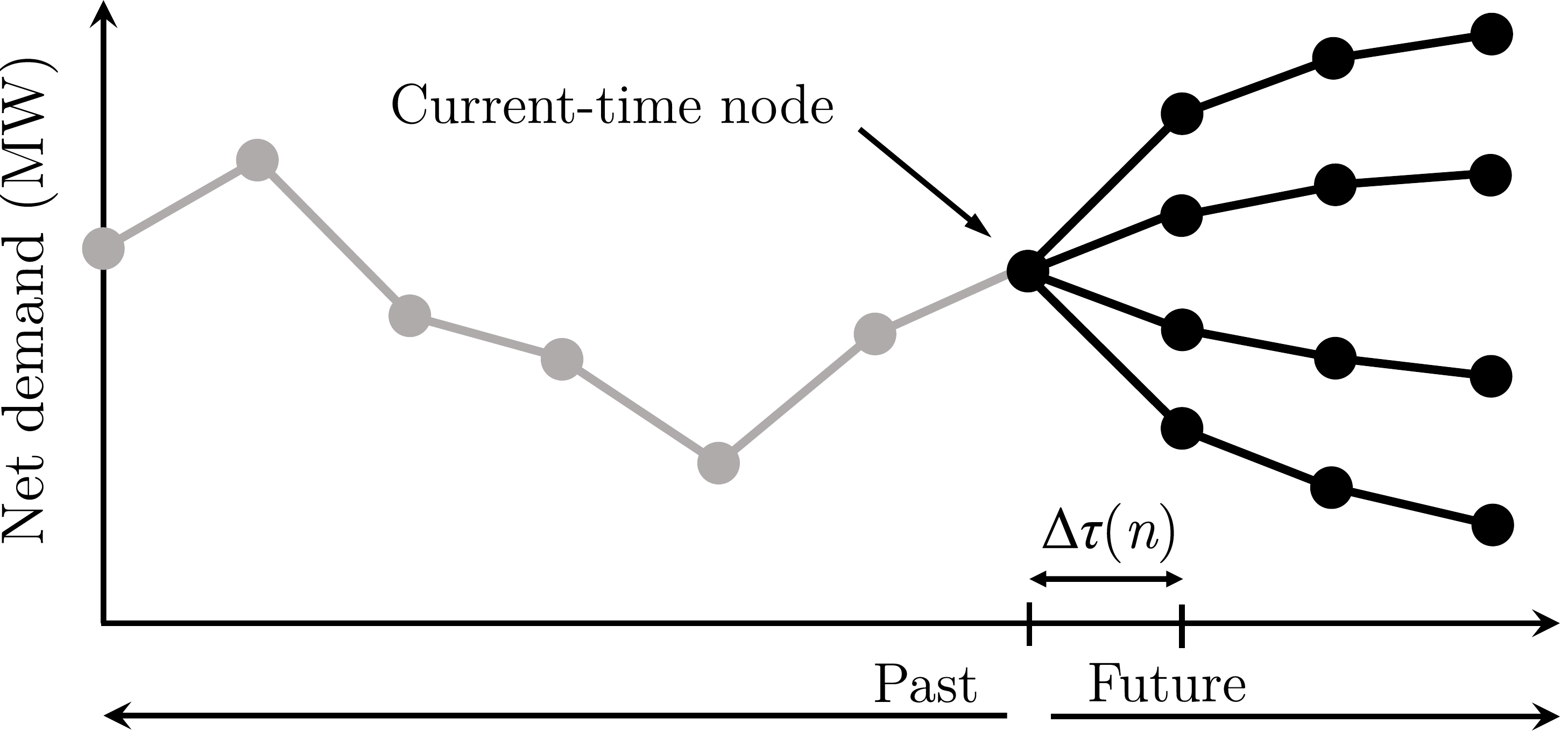}
\caption{Schematic of the scenario tree used in the SUC.}
\label{ScenarioTree}
\end{figure}

A stochastic scheduling model is applied in this paper, which minimises the expected system operation cost while taking into account the uncertainty of RES. Uncertainty is explicitly modelled by means of a scenario tree as that in Fig. \ref{ScenarioTree}, which is built by using the quantile-based scenario generation method described in \cite{AlexEfficient}. Each scenario corresponds to a user-defined quantile of the distribution of the random variable net-demand (demand minus wind power). The study in \cite{AlexEfficient} demonstrated the effectiveness of this approach by selecting a small number of scenarios. From these quantiles, the system conditions and the probability of reaching a particular node in the tree can be derived. For simplicity, trees are constructed with branching only at the current-time node, as this approach has been demonstrated to provide similar results to more intricate tree structures while greatly reducing computational time \cite{AlexEfficient}. 

The Model-Predictive Control concept is applied in this SUC. A full SUC optimisation is firstly computed with a 24h horizon in an hourly time-step $(\Delta\tau(n)=1\textrm{h})$. Only the here-and-now decisions in the current-time node are applied, while all future decisions are discarded. In the next time-step, the realisation of the stochastic variable becomes available, as well as updated forecasts. A new scenario tree is then built, covering another 24-hour horizon, and the SUC optimisation is run again respecting all previously-made decisions and physical constraints. 

The SUC minimises the expected operational cost over all nodes in the scenario tree:
\begin{equation} \label{objectiveSUC}
\text{min}\quad \sum_{n \in \mathcal{N}}\pi(n)\sum_{g \in \mathcal{G}}C_g(n)
\end{equation}
The operating cost of generating units $g$ is given by:
\begin{equation}
C_g(n)=\textrm{c}_g^\textrm{st}N_g^\textrm{sg}(n)+\Delta\tau(n)\left[\textrm{c}_g^\textrm{nl}N_g^\textrm{up}(n)+\textrm{c}_g^\textrm{m}P_g(n)\right]
\end{equation}
Note that generating units with the same characteristics are clustered in the SUC to reduce the computational burden.

The objective function (\ref{objectiveSUC}) is subject to several constraints, in order to correctly model the behaviour of a power system. Frequency-security constraints are discussed in Section \ref{SectionFrequency}. Note that the random variable in the SUC, the net-demand, does not explicitly appear in the mathematical formulation in this paper, but it can be found in the detailed SUC formulation in \cite{AlexEfficient}. Index $(n)$ is dropped from decision variables in following equations presented here to make the expressions clearer.

SUC is a computationally-demanding problem, but it has been shown to provide a more cost-effective operation for systems with high RES penetration. Therefore, our frequency-security model is tested in this demanding problem to demonstrate its applicability. Nevertheless, since it is still common practice for industry to use a Deterministic UC, we also include a comparison of results from Deterministic and Stochastic UC in Section \ref{SectionDetVsStoch}.

\section{Conditions for Complying with Dynamic-Frequency Requirements} \label{SectionFrequency}

Certain dynamic-frequency requirements must be complied with to maintain power system security. In GB, the frequency requirements following a power outage are: 1) RoCoF must be below 0.125Hz/s, 2) the frequency nadir must never be below 49.2Hz and 3) post-fault frequency (also called frequency quasi-steady-state) must recover to be above 49.5Hz after 60s \cite{NationalGridRequirements}. Due to high penetration of nonsynchronous RES, the RoCoF limit is proposed to be relaxed in the near future \cite{NationalGridRelaxRocof}.

In order to maintain frequency security, the scheduling algorithm must be constrained so that a sufficient amount of inertia and FR are available in the system at all times. The deduction of these constraints, as well as some linearisations for them to be implemented in an MILP formulation, are given in the following subsections.

The need for inertia and FR increases along with a higher largest power infeed in the system. The present work considers the option to dynamically reduce the largest power infeed by modelling it as a decision variable in the scheduling problem:
\begin{equation} \label{firstEq}
P_{\textrm{L}} \geq P_g \qquad  \forall \, g \in \mathcal{G}
\end{equation}
This paper uses clustering of equivalent units in the SUC, so the right-hand side of (\ref{firstEq}) must be divided by $N_g^\textrm{up}$. The resulting constraint can be linearised using a big-M method, given that $N_g^\textrm{up}$ is an integer decision variable. For the particular case of the GB system, the largest sources of power are nuclear units, which are modelled as must-run generators. Therefore, the right-hand side of (\ref{firstEq}) becomes $P_\mathcal{M}/\textrm{N}_\mathcal{M}$, yielding a linear constraint.

\subsection{Frequency-Security Constraints}
The constraints that guarantee compliance with dynamic-frequency requirements can be obtained by solving the swing equation, which describes the time-evolution of frequency deviation after a generation outage \cite{KundurBook}:
\begin{multline} \label{SwingEq}
\frac{2H}{f_0}\frac{\textrm{d} \Delta f(t)}{\textrm{d} t}+\textrm{D}\cdot \textrm{P}_{\textrm{D}}\cdot\Delta f(t)\\ 
=\sum_{s \in \mathcal{S}}\textrm{EFR}_s(t)+\sum_{g \in \mathcal{G}}\textrm{PFR}_g(t) - P_{\textrm{L}}
\end{multline}
Where EFR and PFR are modelled as:
\begin{equation} \label{EFRdefinition}
\textrm{EFR}_s(t)=\left\{ 
\begin{array}{ll}
R_s\cdot t / \textrm{T}_{\textrm{s}} \quad & \mbox{if $t\leq\textrm{T}_\textrm{s}$} \\
R_s \quad & \mbox{if $t> \textrm{T}_\textrm{s}$}
\end{array}
\right.
\end{equation}
\begin{equation} \label{PFRdefinition}
\textrm{PFR}_g(t)=\left\{ 
\begin{array}{ll}
R_g\cdot t / \textrm{T}_{\textrm{g}} \quad & \mbox{if $t\leq\textrm{T}_\textrm{g}$} \\
R_g \quad & \mbox{if $t> \textrm{T}_\textrm{g}$}
\end{array}
\right.
\end{equation}

By solving (\ref{SwingEq}), the conditions for respecting dynamic-frequency requirements can be derived. The constraint that guarantees RoCoF security is directly obtained from (\ref{SwingEq}) by realising that the highest value for RoCoF is achieved at the very instant of the outage ($t=0$), when frequency deviation is effectively zero:
\begin{equation} \label{RocofConstraint}
\left|\mbox{RoCoF}\right| = \frac{P_{\textrm{L}}\cdot f_0}{2 H} \leq \mbox{RoCoF}_{\textrm{max}}
\end{equation}
The level of system inertia after the largest loss is given by:
\begin{equation}
H=\sum_{g \in \mathcal{G}}\textrm{H}_g\cdot \textrm{P}_g^{\textrm{max}}\cdot N_g^{\textrm{up}}+\textrm{H}_\textrm{W}\cdot P_{\textrm{W}}- \textrm{P}^{\textrm{max}}_{\textrm{L}}\cdot \textrm{H}_\textrm{L}
\end{equation}

The constraint for assuring quasi-steady-state (q-s-s) security can be obtained from (\ref{SwingEq}) by assuming that RoCoF is effectively zero in quasi-steady-state:
\begin{equation} \label{qssConstraint}
\left| \Delta f^{\textrm{ss}} \right| = \frac{(P_{\textrm{L}} - R_{\mathcal{S}} - R_{\mathcal{G}})}{\textrm{D}\cdot \textrm{P}_{\textrm{D}}} \leq \Delta f^{\textrm{ss}}_{\textrm{max}}
\end{equation}
Note that the total amount of EFR is defined as $R_{\mathcal{S}} = \sum_{s \in \mathcal{S}}R_{s}$ and the total amount of PFR as $R_{\mathcal{G}} = \sum_{g \in \mathcal{G}}R_{g}$.

The nadir requirement is respected if the following inequality holds true:
\begin{equation} \label{nadir_req}
\left| \Delta f_{\textrm{nadir}} \right| = \left| \Delta f(t=t^{*})  \right| \leq \Delta f_{\textrm{max}}
\end{equation}
Where $t^*$ is the time when the nadir is reached. The expression for $\Delta f_{\textrm{nadir}}$ can be obtained by solving (\ref{SwingEq}). However, if damping is considered in (\ref{SwingEq}), it is not possible to obtain an analytical expression for $\Delta f_{\textrm{nadir}}$ as this would imply solving an equation involving sums of exponential functions, which is only solvable by numerical methods. Therefore, here we first obtain the constraint for nadir without considering the effect of load damping, which yields a more conservative condition. Nevertheless, a term approximating the contribution of damping to support the frequency nadir is proposed in Section \ref{SectionDamping}.

In order to obtain $\Delta f_{\textrm{nadir}}$, eq. (\ref{SwingEq}) must be solved for time interval $t \in \left[\textrm{T}_\textrm{s},\textrm{T}_\textrm{g}\right)$: the nadir must certainly occur before $\textrm{T}_\textrm{g}$, as otherwise frequency would drop indefinitely; furthermore, the nadir will take place after $\textrm{T}_s$ if $R_\mathcal{S}$ is lower than the largest power infeed, which is certainly the case for any operating condition in GB's power grid. 
Then, by neglecting the damping term in (\ref{SwingEq}) and solving:
\begin{equation} \label{Deltaf_nadir}
\Delta f(t^{*}) = \frac{f_0}{2 H}\left[\frac{R_{\mathcal{G}}}{2\textrm{T}_g} (t^{*})^2 + R_{\mathcal{S}} \left(t^{*}-\frac{\textrm{T}_s}{2}\right) - P_{\textrm{L}}\cdot t^{*}\right]
\end{equation}
The time at which nadir is reached can be calculated from (\ref{SwingEq}), by setting the derivative of frequency deviation to zero:
\begin{equation} \label{tnadir_simple}
t^* = \frac{(P_{\textrm{L}} - R_{\mathcal{S}})\cdot \textrm{T}_\textrm{g}}{R_{\mathcal{G}}} 
\end{equation}
By substituting (\ref{tnadir_simple}) into (\ref{Deltaf_nadir}), and then substituting the resulting expression into (\ref{nadir_req}), the constraint that guarantees compliance with the nadir requirement can be obtained:
\begin{equation} \label{nadirConstraint}
\left(\frac{H}{f_0} - \frac{R_{\mathcal{S}} \cdot \textrm{T}_\textrm{s}}{4\cdot \Delta f_{\textrm{max}}}\right) \cdot R_{\mathcal{G}} \geq 
\frac{(P_{\textrm{L}}-R_{\mathcal{S}})^2 \cdot \textrm{T}_\textrm{g}}{4\cdot \Delta f_{\textrm{max}}}
\end{equation}

Both FR services considered in this formulation are assumed in (\ref{EFRdefinition}) and (\ref{PFRdefinition}) to start providing response right after the generation outage. Therefore, the frequency deadband of turbine governors and BESS' control systems is not considered in the deduction of the nadir constraint, in order to make following mathematical expressions less intricate. However, if a deadband were to be considered, the procedure presented in \cite{FeiStochastic} could be followed for obtaining the nadir constraint, resulting in an equivalent formulation to (\ref{nadirConstraint}).

The RoCoF constraint (\ref{RocofConstraint}), q-s-s constraint (\ref{qssConstraint}) and nadir constraint (\ref{nadirConstraint}) assure post-fault frequency security in a power grid. By including these constraints in any optimisation routine, such as a UC formulation, a secure post-fault frequency evolution is guaranteed without the need to refine the time-interval for solving the optimisation: the UC can still be solved in the typical hourly or half-hourly basis. The frequency-security constraints enforce that a sufficient amount of resources such as inertia and FR are scheduled so that, if a power outage were to occur during this time interval, the sub-second dynamics of frequency are guaranteed to be within limits.

\subsection{Contribution of Load Damping to Supporting the Nadir} \label{SectionDamping}

\begin{figure}[!tp]
\centering
%
%
\definecolor{mycolor1}{rgb}{1.00000,0.26275,0.26275}%
\begin{tikzpicture}

\begin{axis}[%
axis lines = left, 
width=3in,
height=1.05in,
at={(0in,0in)},
scale only axis,
clip=false,
xmin=0,
xmax=117.5+7,
xtick={0,112.5},
xticklabels={{},{}},
ymin=0,
ymax=10.925+0.5,
ytick={0,10.125},
yticklabels={{},{}},
axis background/.style={fill=white},
title style={font=\Huge},xlabel style={font={\color{blue}\bfseries}},ylabel style={font=\tiny},legend style={font=\scriptsize},ticklabel style={font=\color{red}}
]
\addplot [color=black, line width=1.0pt]
  table[row sep=crcr]{%
1e-06	10.1250000539445\\
1.000001	10.0051780807878\\
2.000001	9.88571438316904\\
3.000001	9.76661136129624\\
4.000001	9.64787162059976\\
5.000001	9.52949779847449\\
6.000001	9.41149256498808\\
7.000001	9.29385862361181\\
8.000001	9.17659871197638\\
9.000001	9.05971560265215\\
10.000001	8.9432121039559\\
11.000001	8.82709106078494\\
12.000001	8.71135535547986\\
13.000001	8.59600790871699\\
14.000001	8.4810516804322\\
15.000001	8.36648967077731\\
16.000001	8.25232492111076\\
17.000001	8.13856051502382\\
18.000001	8.02519957940442\\
19.000001	7.91224528554015\\
20.000001	7.79970085026237\\
21.000001	7.68756953713344\\
22.000001	7.5758546576792\\
23.000001	7.46455957266894\\
24.000001	7.35368769344529\\
25.000001	7.24324248330649\\
26.000001	7.13322745894392\\
27.000001	7.02364619193748\\
28.000001	6.9145023103122\\
29.000001	6.80579950015913\\
30.000001	6.69754150732405\\
31.000001	6.58973213916772\\
32.000001	6.48237526640168\\
33.000001	6.37547482500375\\
34.000001	6.26903481821779\\
35.000001	6.16305931864261\\
36.000001	6.05755247041525\\
37.000001	5.95251849149405\\
38.000001	5.84796167604773\\
39.000001	5.74388639695666\\
40.000001	5.64029710843349\\
41.000001	5.53719834877036\\
42.000001	5.4345947432208\\
43.000001	5.33249100702499\\
44.000001	5.23089194858773\\
45.000001	5.12980247281911\\
46.000001	5.02922758464891\\
47.000001	4.92917239272654\\
48.000001	4.8296421133193\\
49.000001	4.73064207442298\\
50.000001	4.63217772009993\\
51.000001	4.53425461506104\\
52.000001	4.43687844950972\\
53.000001	4.34005504426747\\
54.000001	4.24379035620255\\
55.000001	4.14809048398515\\
56.000001	4.05296167419504\\
57.000001	3.95841032780985\\
58.000001	3.86444300710524\\
59.000001	3.77106644300115\\
60.000001	3.67828754289223\\
61.000001	3.5861133990042\\
62.000001	3.49455129732258\\
63.000001	3.40360872714544\\
64.000001	3.3132933913173\\
65.000001	3.22361321720836\\
66.000001	3.13457636851016\\
67.000001	3.04619125792781\\
68.000001	2.9584665608584\\
69.000001	2.87141123015647\\
70.000001	2.78503451210052\\
71.000001	2.69934596368914\\
72.000001	2.61435547141289\\
73.000001	2.53007327166805\\
74.000001	2.44650997300163\\
75.000001	2.36367658040478\\
76.000001	2.28158452190368\\
77.000001	2.20024567773523\\
78.000001	2.11967241243986\\
79.000001	2.03987761025755\\
80.000001	1.96087471427758\\
81.000001	1.88267776986938\\
82.000001	1.80530147301594\\
83.000001	1.7287612242844\\
84.000001	1.65307318930769\\
85.000001	1.57825436682233\\
86.000001	1.50432266551947\\
87.000001	1.43129699123173\\
88.000001	1.35919734631132\\
89.000001	1.28804494347885\\
90.000001	1.21786233696394\\
91.000001	1.14867357446007\\
92.000001	1.08050437433226\\
93.000001	1.01338233372769\\
94.000001	0.947337174860279\\
95.000001	0.88240103894112\\
96.000001	0.818608840256411\\
97.000001	0.755998697139368\\
98.000001	0.694612462636277\\
99.000001	0.634496386481669\\
100.000001	0.575701953131838\\
101.000001	0.518286960694505\\
102.000001	0.462316937250791\\
103.000001	0.407867042709322\\
104.000001	0.355024692048348\\
105.000001	0.303893292092177\\
106.000001	0.25459777901393\\
107.000001	0.207293241948985\\
108.000001	0.162179246098655\\
109.000001	0.119525793039578\\
110.000001	0.0797267074238576\\
111.000001	0.0434342011273137\\
112.000001	0.0120737923622273\\
};

\addplot [color=black, dashed, line width=1.0pt]
  table[row sep=crcr]{%
1e-06	10.12499991\\
1.000001	10.03499991\\
2.000001	9.94499991\\
3.000001	9.85499991\\
4.000001	9.76499991\\
5.000001	9.67499991\\
6.000001	9.58499991\\
7.000001	9.49499991\\
8.000001	9.40499991\\
9.000001	9.31499991\\
10.000001	9.22499991\\
11.000001	9.13499991\\
12.000001	9.04499991\\
13.000001	8.95499991\\
14.000001	8.86499991\\
15.000001	8.77499991\\
16.000001	8.68499991\\
17.000001	8.59499991\\
18.000001	8.50499991\\
19.000001	8.41499991\\
20.000001	8.32499991\\
21.000001	8.23499991\\
22.000001	8.14499991\\
23.000001	8.05499991\\
24.000001	7.96499991\\
25.000001	7.87499991\\
26.000001	7.78499991\\
27.000001	7.69499991\\
28.000001	7.60499991\\
29.000001	7.51499991\\
30.000001	7.42499991\\
31.000001	7.33499991\\
32.000001	7.24499991\\
33.000001	7.15499991\\
34.000001	7.06499991\\
35.000001	6.97499991\\
36.000001	6.88499991\\
37.000001	6.79499991\\
38.000001	6.70499991\\
39.000001	6.61499991\\
40.000001	6.52499991\\
41.000001	6.43499991\\
42.000001	6.34499991\\
43.000001	6.25499991\\
44.000001	6.16499991\\
45.000001	6.07499991\\
46.000001	5.98499991\\
47.000001	5.89499991\\
48.000001	5.80499991\\
49.000001	5.71499991\\
50.000001	5.62499991\\
51.000001	5.53499991\\
52.000001	5.44499991\\
53.000001	5.35499991\\
54.000001	5.26499991\\
55.000001	5.17499991\\
56.000001	5.08499991\\
57.000001	4.99499991\\
58.000001	4.90499991\\
59.000001	4.81499991\\
60.000001	4.72499991\\
61.000001	4.63499991\\
62.000001	4.54499991\\
63.000001	4.45499991\\
64.000001	4.36499991\\
65.000001	4.27499991\\
66.000001	4.18499991\\
67.000001	4.09499991\\
68.000001	4.00499991\\
69.000001	3.91499991\\
70.000001	3.82499991\\
71.000001	3.73499991\\
72.000001	3.64499991\\
73.000001	3.55499991\\
74.000001	3.46499991\\
75.000001	3.37499991\\
76.000001	3.28499991\\
77.000001	3.19499991\\
78.000001	3.10499991\\
79.000001	3.01499991\\
80.000001	2.92499991\\
81.000001	2.83499991\\
82.000001	2.74499991\\
83.000001	2.65499991\\
84.000001	2.56499991\\
85.000001	2.47499991\\
86.000001	2.38499991\\
87.000001	2.29499991\\
88.000001	2.20499991\\
89.000001	2.11499991\\
90.000001	2.02499991\\
91.000001	1.93499991\\
92.000001	1.84499991\\
93.000001	1.75499991\\
94.000001	1.66499991\\
95.000001	1.57499991\\
96.000001	1.48499991\\
97.000001	1.39499991\\
98.000001	1.30499991\\
99.000001	1.21499991\\
100.000001	1.12499991\\
101.000001	1.03499991\\
102.000001	0.94499991\\
103.000001	0.85499991\\
104.000001	0.76499991\\
105.000001	0.67499991\\
106.000001	0.58499991\\
107.000001	0.49499991\\
108.000001	0.40499991\\
109.000001	0.31499991\\
110.000001	0.22499991\\
111.000001	0.13499991\\
112.000001	0.0449999100000002\\
};

\node[right, align=left]
at (axis cs:-3,-1.2) {\footnotesize 0};
\node[right, align=left]
at (axis cs:100,7) {$\frac{(P_{\textrm{L}}-R_\mathcal{S})}{\textrm{D} \cdot \Delta f_{\textrm{max}}}$};
\node[right, align=left] (C)
at (axis cs:110.5,6) {};
\node[right, align=left] (D)
at (axis cs:110.5,0.5) {};
\draw [->,>=stealth] (C) -- (D);
\node[right, align=left]
at (axis cs:-8,0.1) {\footnotesize 0};
\node[right, align=left] (A)
at (axis cs:30,10.225) {$\frac{(P_{\textrm{L}}-R_\mathcal{S})^2 \cdot\textrm{T}_\textrm{g}\cdot f_0}{4 H\cdot\Delta f_{\textrm{max}}-R_\mathcal{S}\cdot \textrm{T}_\textrm{s}\cdot f_0}$};
\node[right, align=left] (B)
at (axis cs:5,10.225) {};
\draw [->,>=stealth] (A) -- (B);
\node[right, align=left, rotate=90]
at (axis cs:-8,2.5) {\small $R_\mathcal{G}  \; (\mbox{MW})$};
\node[right, align=left]
at (axis cs:46.25,-1.8) {\small $\textrm{P}_\textrm{D} \; (\mbox{MW})$};

\end{axis}
\end{tikzpicture}%
\caption{Exact feasible region for respecting the nadir requirement, given by the epigraph of the solid curve; and feasible region defined by the proposed linear approximation, given by the epigraph of the dashed line.}
\label{FigNadirLinear}
\end{figure}
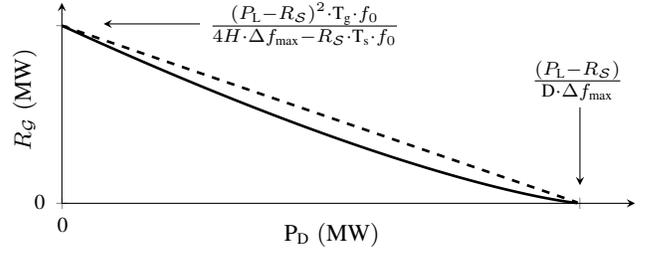

Here we propose a linear term to be included in (\ref{nadirConstraint}) in order to account for the effect of load damping. First, $\Delta f(t^*)$ is obtained while considering damping in (\ref{SwingEq}):
\begin{multline} \label{NadirExact}
\Delta f(t^*)\cdot\frac{\mathrm{D}'}{f_0}=\left(\textrm{e}^{-\frac{\mathrm{D}'}{2H} t^*}-1\right)\left(P_{\textrm{L}}+\frac{2H\cdot R_\mathcal{G}}{\mathrm{D}'\cdot\textrm{T}_\textrm{g}}\right)+\frac{R_\mathcal{G}}{\textrm{T}_\textrm{g}} t^* \\ 
+ R_\mathcal{S}\left[1+\frac{2H}{\mathrm{D}'\cdot\textrm{T}_\textrm{s}}\left(\textrm{e}^{-\frac{\mathrm{D}'}{2H} t^*}-\textrm{e}^{-\frac{\mathrm{D}'}{2H} (t^*-\textrm{T}_\textrm{s})}\right)\right]
\end{multline}
Where $\mathrm{D}'=\mathrm{D}\cdot \mathrm{P}_\textrm{D} \cdot f_0$ and $t^*$ can be obtained from setting RoCoF to zero in (\ref{SwingEq}):
\begin{equation} \label{tnadir}
t^* = \frac{(P_{\textrm{L}} - R_{\mathcal{S}} - \mathrm{D}\cdot \mathrm{P}_\textrm{D}\cdot \left|\Delta f_\textrm{nadir}\right|)\cdot \textrm{T}_\textrm{g}}{R_{\mathcal{G}}} 
\end{equation}

Since (\ref{NadirExact}) implicitly describes a monotonically increasing function $\Delta f(R_\mathcal{G})$,
the necessary $R_\mathcal{G}$ for a particular system condition can be obtained by numerically solving the following equation:
\begin{multline} \label{NadirSolveNumerical}
-\Delta f_\textrm{max}\cdot\frac{\mathrm{D}'}{f_0}=\left(\textrm{e}^{-\gamma}-1\right)\left(P_{\textrm{L}}+\frac{2H\cdot R_\mathcal{G}}{\mathrm{D}'\cdot\textrm{T}_\textrm{g}}\right)+P_{\textrm{L}} - R_{\mathcal{S}} 
\\ 
- \frac{\mathrm{D}'}{f_0}\cdot \Delta f_\textrm{max}
+ R_\mathcal{S}\left[1+\frac{2H}{\mathrm{D}'\cdot\textrm{T}_\textrm{s}}\left(\textrm{e}^{-\gamma}-\textrm{e}^{-\gamma+\frac{\mathrm{D}'}{2H}\textrm{T}_\textrm{s}}\right)\right]
\end{multline}
Where $\gamma$ is defined as:
\begin{equation} 
\gamma = \frac{\mathrm{D}\cdot \mathrm{P}_\textrm{D}\cdot f_0 \cdot \textrm{T}_\textrm{g} \cdot (P_{\textrm{L}} - R_{\mathcal{S}} - \mathrm{D}\cdot \mathrm{P}_\textrm{D}\cdot \Delta f_\textrm{max})}{2H\cdot R_{\mathcal{G}}}
\end{equation}

Eq. (\ref{NadirSolveNumerical}) implicitly describes a convex function $R_{\mathcal{G}}(\mathrm{P}_\textrm{D})$, as shown by the solid curve in Fig. \ref{FigNadirLinear}, which represents the numerical solution of (\ref{NadirSolveNumerical}) for several values of $\mathrm{P}_\textrm{D}$. As $R_{\mathcal{G}}$ and $\textrm{P}_{\textrm{D}}$ are both positive physical magnitudes, the range of $\textrm{P}_{\textrm{D}}$ for which (\ref{NadirSolveNumerical}) has a solution is: 
\begin{equation}
\textrm{P}_{\textrm{D}} \in \left(0, \; \frac{(P_{\textrm{L}}-R_\mathcal{S})}{\textrm{D}\cdot \Delta f_{\textrm{max}}}\right)
\end{equation}

As the function $R_{\mathcal{G}}(\mathrm{P}_\textrm{D})$ described by (\ref{NadirSolveNumerical}) is convex and monotonically decreasing, it can be inner-approximated by a line. 
A graphical representation of this linear approximation is also provided in Fig. \ref{FigNadirLinear}. Therefore, the linearised effect of load damping on nadir can be included in the previously obtained nadir constraint (\ref{nadirConstraint}) as follows:
\begin{equation} \label{nadirEFRDamping}
\left(\frac{H}{f_0} - \frac{R_{\mathcal{S}} \cdot \textrm{T}_\textrm{s}}{4\cdot \Delta f_{\textrm{max}}}\right) \cdot R_{\mathcal{G}} \geq 
\alpha -\beta \cdot \textrm{P}_{\textrm{D}}
\end{equation}
Where:
\begin{equation} \label{alphaEq}
\alpha = \frac{(P_{\textrm{L}}-R_{\mathcal{S}})^2 \cdot \textrm{T}_\textrm{g}}{4 \cdot \Delta f_{\textrm{max}}}, \quad
\beta = \frac{(P_{\textrm{L}}-R_{\mathcal{S}})\cdot \textrm{T}_\textrm{g} \cdot \textrm{D}}{4}
\end{equation}

The inner approximation leads to an underestimation of the actual contribution from load damping, as can be clearly observed in Fig. \ref{FigNadirLinear}: the feasible region defined by the linear approximation, which is the epigraph of the dashed line, is tighter than the actual feasible region, epigraph of the solid curve. Nevertheless, this underestimation is significantly less conservative than simply neglecting the effect of damping. The quantitative assessment of this proposed linear approximation is presented in Section \ref{QuantitativeApproxSection}.

\subsection{Linearisation of the Frequency-Nadir Constraint} \label{sectionLinearisation}

\begin{figure}[t!] 
  \centering
  \includegraphics[width=3.25in]{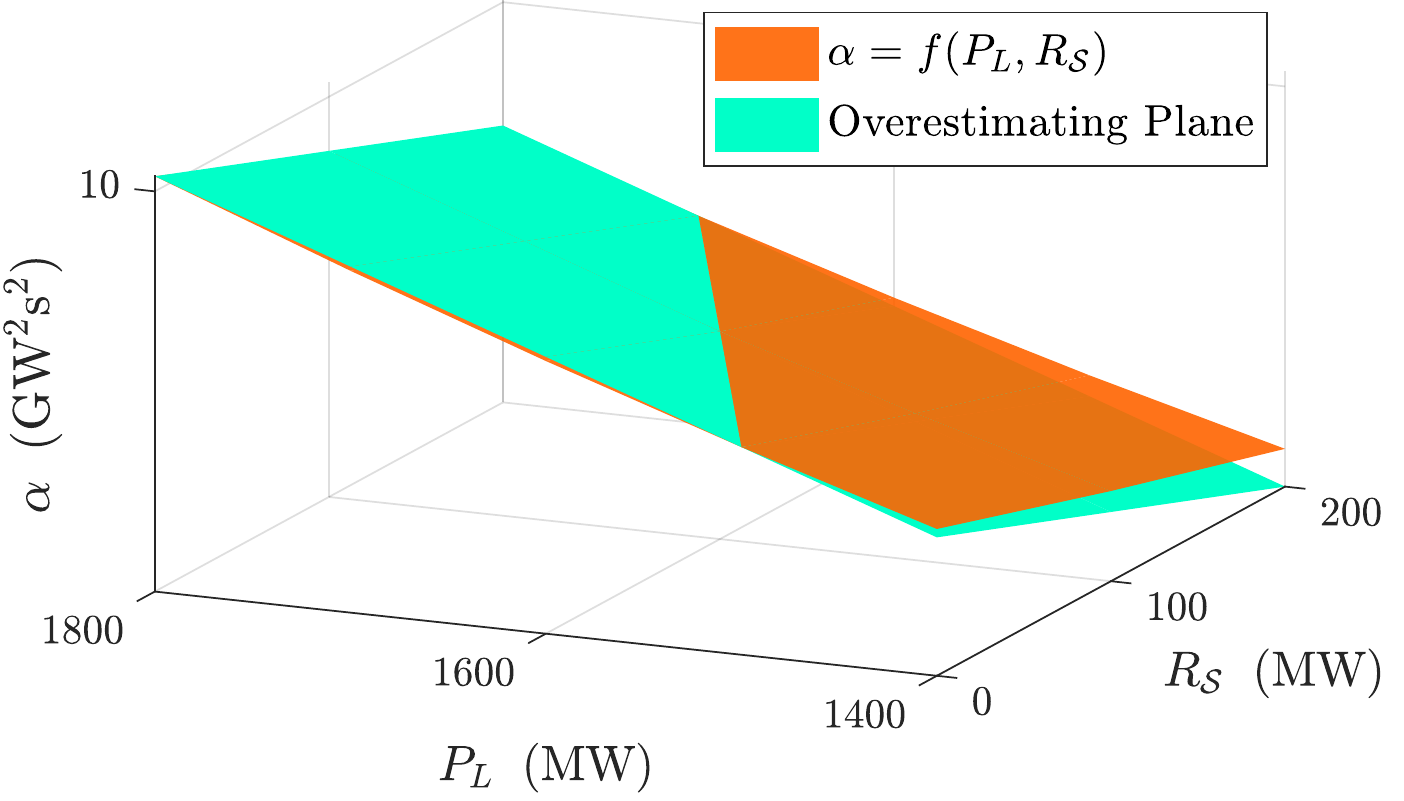}
  \caption{Representation of the surface defined by $\alpha$ in (\ref{alphaEq}) and a plane which overestimates the surface in one-half of its domain.}
  \label{FigPlanes}
\end{figure}

\begin{figure*}[!bh]
\noindent\rule{\textwidth}{0.7pt}
\vspace{4pt}
    \begin{equation} \label{nadirPlanes}
\left(\frac{H}{f_0} - \frac{R_{\mathcal{S}} \cdot \textrm{T}_\textrm{s}}{4\cdot \Delta f_{\textrm{max}}}\right) \cdot R_{\mathcal{G}} \,\, \geq \,\, \textrm{a}_p \cdot P_{\textrm{L}} + \textrm{b}_p \cdot R_{\mathcal{S}} + \textrm{c}_p  
\, - \frac{(P_{\textrm{L}}-R_{\mathcal{S}})\cdot \textrm{T}_\textrm{g} \cdot \textrm{D}}{4} \cdot \textrm{P}_{\textrm{D}}  \qquad \forall \, p \in \mathcal{P}
    \end{equation}
    \vspace*{-12pt}
\end{figure*}

\stepcounter{equation}
\begin{figure*}[!bp] 
    \begin{equation} \label{nadirFINAL}
\frac{1}{f_0}\sum_{l \in \mathcal{L}}{m_l \cdot 2^l} - \frac{\textrm{T}_\textrm{s}}{4\cdot \Delta f_{\textrm{max}}} \cdot \sum_{l \in \mathcal{L}} k_l \cdot 2^l \\
\,\, \geq \,\, \textrm{a}_p \cdot P_{\textrm{L}} + \textrm{b}_p \cdot R_{\mathcal{S}} + \textrm{c}_p  
\, - \frac{(P_{\textrm{L}}-R_{\mathcal{S}})\cdot \textrm{T}_\textrm{g} \cdot \textrm{D}}{4} \cdot \textrm{P}_{\textrm{D}}  \qquad \forall \, p \in \mathcal{P}
    \end{equation}
\end{figure*}

The RoCoF and q-s-s constraints, (\ref{RocofConstraint}) and (\ref{qssConstraint}) respectively, are linear. However, the nadir constraint with the proposed damping term (\ref{nadirEFRDamping}) is nonconvex. An efficient linearisation method is proposed so that this constraint can be implemented in an MILP while guaranteeing frequency security in all cases.

Constraint (\ref{nadirEFRDamping}) contains three nonlinear terms: two products of continuous variables on its left-hand side, namely $H\cdot R_{\mathcal{G}}$ and $R_{\mathcal{S}} \cdot R_{\mathcal{G}}$, and the quadratic term $(P_{\textrm{L}}-R_{\mathcal{S}})^2$ on its right-hand side. Although it is not possible to exactly linearise these terms, we propose a method to approximate them by linear expressions to any desired level of accuracy. This method consists of the following steps: 1) for the right-hand side of the constraint, use an inner approximation by overestimating planes as in \cite{InnerApproxPaper}, applicable as the quadratic term is a convex function; 2) for the left-hand side, the technique proposed in \cite{binaryApprox} can be applied. This technique consists on representing one of the continuous variables by its binary expansion, so that the product of two continuous variables is converted into products of a continuous and several binary variables, which can then be exactly linearised. In this work, the decision variable chosen to be represented by its binary expansion is $R_{\mathcal{G}}$, since it appears in both products on the left-hand side of constraint (\ref{nadirEFRDamping}). 

Using the inner approximation with overestimating planes for the right-hand side of (\ref{nadirEFRDamping}), (\ref{nadirPlanes}) is obtained. A graphical example using two planes is presented in Fig. \ref{FigPlanes}, where each of the planes covers one-half of the domain of the squared function $\alpha$ defined in (\ref{alphaEq}). Only one of the two planes has been graphed for clarity.

    
Regarding the linearisation of the left-hand side of (\ref{nadirEFRDamping}), the binary expansion of $R_\mathcal{G}$ is defined as:
\addtocounter{equation}{-2}
\begin{equation} \label{binaryExpansion}
R_{\mathcal{G}} = \sum_{l \in \mathcal{L}} z_l \cdot 2^l
\end{equation}
After using this binary expansion, the product of continuous variables $H$ and $R_{\mathcal{S}}$ by each of the binary variables $z_l$ can be linearised using the standard big-M technique described in \cite{bigMbook}, which yields the fully linearised nadir constraint (\ref{nadirFINAL}).
%
%
In (\ref{nadirFINAL}), $m_l=H \cdot z_l$ and $k_l=R_{\mathcal{S}} \cdot z_l$, and the auxiliary variables $m_l$ and $k_l$ must be appropriately constrained following the big-M method. 


Both of the approximations proposed in this subsection can be scaled up to achieve any level of accuracy desired, with the tradeoff of a higher computational burden when solving the scheduling problem. The inner approximation by overestimating planes becomes more precise with increasing number of planes, while it adds a higher number of constraints to the optimisation. The binary expansion of $R_{\mathcal{G}}$
provides a better approximation of the continuous variable when a higher number of binary decision variables is used. 
A quantification of the tradeoff between accuracy and computational burden is presented in Section \ref{QuantitativeApproxSection}.

\section{Case Studies} \label{SectionCaseStudies}

\begin{table}[!t]
\renewcommand{\arraystretch}{1.05}
\caption{Characteristics of Thermal Plants in GB's 2030 System}
\label{TableThermal}
\centering
\begin{tabular}{l| l l l l}
    \multicolumn{1}{c|}{} & Nuclear & CCGT & OCGT\\
\hline
Number of Units & 4 & 100 & 30\\
Rated Power (MW) & 1800 & 500 & 100\\
Min Stable Generation (MW) & 1400 & 250 & 50\\
No-Load Cost $\textrm{c}^{\textrm{nl}}_g$ (\pounds/h) & 0 & 4500 & 3000\\
Marginal Cost $\textrm{c}^{\textrm{m}}_g$ (\pounds/MWh) & 10 & 47 & 200\\
Startup Cost $\textrm{c}^{\textrm{st}}_g$ (\pounds) & N/A & 10000 & 0\\
Startup Time (h) & N/A & 4 & 0\\
Min Up Time (h) & N/A & 4 & 0\\
Min Down Time (h) & N/A & 1 & 0\\
Inertia Constant $\textrm{H}_g$ (s) & 5 & 4 & 4\\
Max $R_g$ deliverable (MW) & 0 & 50 & 20\\
Emissions (ton$\mbox{CO}_{2}$/MWh) & 0 & 0.394 & 0.557\\
\hline
\end{tabular}
\end{table}

The proposed model can be used to achieve cost-effective operation of a low-carbon power system. In addition, it would also allow to identify the technologies and practices that would be most beneficial to the system, such as BESS providing EFR, or part-loading the largest generating units. The benefits from operating a system using this frequency-secured framework are not only economic but also environmental, as shown through case studies in this section.

Several case studies on the GB 2030 power system are carried out to demonstrate the benefits of the proposed model. The characteristics of generation plants are included in Table \ref{TableThermal}. Generation costs are based on projected $\textrm{CO}_2$ cost of \pounds45/ton in 2030. 
The minimum and maximum demand are 20GW and 60GW, respectively. A pumped storage unit with 10GWh capacity, 2.6GW rating and 75\% round efficiency is also present in the system. Moreover, BESS with 1GWh capacity, 200MW rating and 90\% efficiency is installed with the capability to provide EFR. Other system's parameters are: $\textrm{T}_\textrm{s}=0.5\textrm{s}$, $\textrm{T}_\textrm{g}=10\textrm{s}$, $\textrm{D}=0.5\textrm{\%/Hz}$ and $\textrm{P}^{\textrm{max}}_{\textrm{L}}=1800\textrm{MW}$. The dynamic-frequency requirements are set by National Grid regulation to $\Delta f_\textrm{max}=0.8\textrm{Hz}$, $\Delta f^\textrm{ss}_\textrm{max}=0.5\textrm{Hz}$ and $\textrm{RoCoF}_\textrm{max}=0.5\textrm{Hz/s}$ (corresponding to the expected relaxed $\textrm{RoCoF}$ requirement). Synthetic inertia from wind turbines is only considered in Section \ref{sectionSI}.

Simulations were run in a twelve-core 3.5GHz Intel Xeon CPU with 64GB of RAM. 
The optimisations were solved with FICO Xpress 8.0, linked to a multi-threaded C++ application via the BCL interface. The duality gap for the MILPs was set to 0.1\%. For the SUC, a scenario tree branching only in the current-time node was used, with net-demand quantiles of 0.005, 0.1, 0.3, 0.5, 0.7, 0.9 and 0.995. This approach using few quantiles was proven by \cite{AlexEfficient} to provide similar results to more intricate structures, while significantly reducing the computational burden. The values for the quantiles were chosen to capture symmetrically the variance of the forecast error distribution, as well as its tails. Each SUC simulation corresponds to four months, in different seasons, of the operation of GB's 2030 system.


\subsection{Validation of the Frequency-Security Constraints}

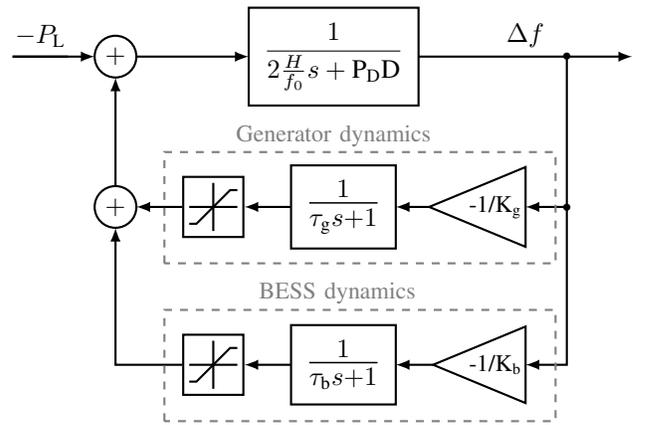
\begin{figure} 
\raggedright

\tikzset{%
  block/.style    = {draw, thick, rectangle, minimum height = 3em,
    minimum width = 3em, inner sep=2mm, outer sep=0mm,},
  gain/.style     = {draw, thick, isosceles triangle, minimum height = 0.5em, inner sep=0.5mm, outer sep=0mm,
     isosceles triangle apex angle=45, shape border rotate=-180},
  sum/.style      = {draw, circle, inner sep=2pt, outer sep=0pt, node distance = 1cm}, 
  input/.style    = {coordinate}, 
  output/.style   = {coordinate}, 
  saturation/.style={%
    draw, 
    path picture={
      \pgfpointdiff{\pgfpointanchor{path picture bounding box}{north east}}%
        {\pgfpointanchor{path picture bounding box}{south west}}
      \pgfgetlastxy\x\y
      \tikzset{x=\x*.4, y=\y*.4}
      %
      \draw (-1,0) -- (1,0) (0,-1) -- (0,1); 
      \draw (-1,-.7) -- (-.7,-.7) -- (.7,.7) -- (1,.7);
    }
  }
}
\newcommand{\suma}{$+$}

\begin{tikzpicture}[auto, thick, node distance=1cm, >=latex]

\draw
    node [] (input1) {}
    node [right of=input1] (write_PLoss) {}
    node at (1.2,0)[] (before_sum1) {}
    node at (1.5,0)[sum] (sum1) {\suma}
    node at (4.4,0)[block] (InertiaBlock) {$\dfrac{1}{2\frac{H}{f_0} s + \textrm{P}_\textrm{D} \textrm{D} }$}
    node at (8.5,0)[] (output1) {}
    node at (7.5,0.1)[](derivation1) {}
    node at (6.55,-2)[gain] (droop1) {\small -1/$\textrm{K}_\textrm{g}$}
    node at (4.53,-2)[block] (governor) {\Large $\frac{1}{\tau_\textrm{g} s + 1}$}
    node at (2.8,-2)[saturation, minimum width=0.8cm, minimum height=0.8cm] (sat1) {}
    node at (7.5,-1.9)[](derivation2) {}
    node at (6.55,-4.1)[gain] (droop2) {\small -1/$\textrm{K}_\textrm{b}$}
    node at (4.53,-4.1)[block] (battery) {\Large $\frac{1}{\tau_\textrm{b} s + 1}$}
    node at (2.8,-4.1)[saturation, minimum width=0.8cm, minimum height=0.8cm] (sat2) {}
    node at (1.5,-2)[sum](sum2) {\suma}
    ;
    
	\draw[](input1) -- node {$-P_\textrm{L}$} (write_PLoss);
    \draw[](input1) -- node {} (before_sum1);
    \draw[->](write_PLoss) -- node {} (sum1);
	\draw[->](sum1) -- node {} (InertiaBlock);
    \draw[->](InertiaBlock) -- node {$\Delta f$} (output1);
    \draw[->](derivation1) |- node[near end]{} (droop1);
    \draw[->](droop1) -- node[near end]{} (governor);
	\draw[->](governor) -- node[near end]{} (sat1);
    \draw[->](sat1) -- node[near end]{} (sum2);
    \draw[->](derivation2) |- node[near end]{} (droop2);
    \draw[->](droop2) -- node[near end]{} (battery);
	\draw[->](battery) -- node[near end]{} (sat2);
    \draw[->](sat2) -| node[near end]{} (sum2);
    \draw[->](sum2) -- node[near end]{} (sum1);

\draw
	node at (7.5,-0.035){\textbullet}
    node at (7.5,-2.035){\textbullet}
    ;
    
	\draw [color=gray, dashed,thick](2.15,-2.75) rectangle (7.35,-1.27);
	\node [color=gray] at (4.4,-1.05) {\small Generator dynamics};
    \draw [color=gray, dashed,thick](2.15,-4.85) rectangle (7.35,-3.37);
	\node [color=gray] at (4.45,-3.15) {\small BESS dynamics};
	
\end{tikzpicture}
\caption{Block diagram for the simulation of the system frequency dynamics.}
\label{FigBlockDiagram}
\end{figure}

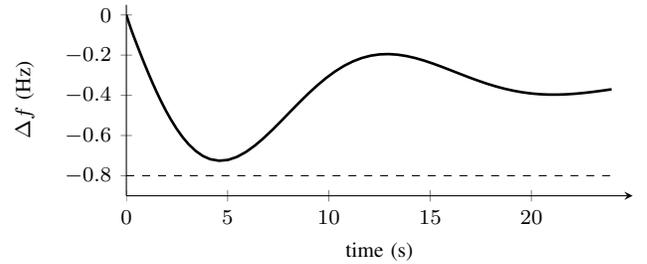
\begin{figure}
\hspace*{0.05cm}
%
%
\definecolor{mycolor1}{rgb}{1.00000,0.26275,0.26275}%
\begin{tikzpicture}

\begin{axis}[%
axis lines = left,
width=2.65in,
height=1in,
at={(0in,0in)},
scale only axis,
y label style={at={(axis description cs:-0.2,0.25)},anchor=west},
xmin=0,
xmax=25,
xlabel style={font=\footnotesize},
xlabel={time (s)},
xtick={0,5,10,15,20},
ymin=-0.9,
ymax=0.05,
ylabel style={font=\footnotesize},
ylabel style={align=center},
ylabel={$\Delta f$ (Hz)},
axis background/.style={fill=white},
legend style={at={(0.03,0.97)}, anchor=north west, legend cell align=left, align=left, draw=white!15!black},legend style={font=\footnotesize},
legend style={fill opacity=0,text opacity=1,draw=none}, 
xticklabel style={font=\footnotesize},
yticklabel style={font=\footnotesize},
y axis line style = {-} 
]






\addplot [color=black, line width=1.0pt]
  table[row sep=crcr]{%
0	0
0.000878510611805197	-0.000276188391203802\\
0.00263553183541559	-0.000828542710816647\\
0.00439255305902598	-0.00138084731035192\\
0.00659677917761217	-0.00207364149589513\\
0.0141905719478203	-0.00445937184641769\\
0.0217843647180285	-0.00684307099570473\\
0.0293781574882366	-0.00922407172206726\\
0.0369719502584447	-0.0116017552287181\\
0.0445657430286529	-0.0139755477766463\\
0.0521595357988610	-0.0163449175756708\\
0.0597533285690692	-0.0187093718987255\\
0.0695872156156866	-0.0217632429293480\\
0.0794211026623040	-0.0248072066048925\\
0.0919020512817114	-0.0286549455869056\\
0.104382999901119	-0.0324839221880852\\
0.119979708414718	-0.0372403468474873\\
0.135576416928317	-0.0419632309578865\\
0.155159700144437	-0.0478428679366945\\
0.174742983360557	-0.0536636619258564\\
0.199619061374250	-0.0609688917201289\\
0.224495139387942	-0.0681714448195172\\
0.256699985228390	-0.0773388331804033\\
0.288904831068839	-0.0863254803722722\\
0.331905154101315	-0.0980388349353224\\
0.344354639292196	-0.101368830966056\\
0.344354639292200	-0.101368830966057\\
0.353347890693980	-0.103762368370300\\
0.371334393497540	-0.108545329226904\\
0.389320896301100	-0.113320047095175\\
0.412604959739001	-0.119488442790916\\
0.521967357756700	-0.148252985340754\\
0.631329755774400	-0.176646396446787\\
0.784402072933978	-0.215679290882472\\
0.937474390093556	-0.253792463334916\\
1.12157852769534	-0.298252435824397\\
1.30568266529711	-0.341054236775588\\
1.48978680289889	-0.382044354677357\\
1.72453490392739	-0.431434241853709\\
1.95928300495589	-0.477381772132361\\
2.19403110598439	-0.519688276497778\\
2.48598189783308	-0.566972425417105\\
2.77793268968177	-0.608146955351297\\
3.06988348153045	-0.643076556698259\\
3.42725262548278	-0.677280406712333\\
3.78462176943511	-0.702161573043407\\
4.14199091338744	-0.717970529867376\\
4.58235309604215	-0.725679042123244\\
5.02271527869685	-0.721391124659455\\
5.46307746135155	-0.706387006240173\\
5.90343964400625	-0.682158547043692\\
6.34380182666095	-0.650321199759509\\
6.78416400931565	-0.612554570520387\\
7.22452619197035	-0.570546969039101\\
7.66488837462506	-0.525944918207288\\
8.10525055727976	-0.480309260197344\\
8.54561273993446	-0.435077784893230\\
8.98597492258916	-0.391535247572420\\
9.42633710524387	-0.350790766990379\\
9.86669928789857	-0.313762601969723\\
10.3070614705533	-0.281170007140306\\
10.7474236532080	-0.253531716936461\\
11.1877858358627	-0.231170449787195\\
11.6281480185174	-0.214222709607939\\
12.0685102011721	-0.202653077338994\\
12.5088723838268	-0.196272133038036\\
12.9492345664815	-0.194757127196243\\
13.3895967491362	-0.197674525971222\\
13.8299589317909	-0.204503586138892\\
14.2703211144456	-0.214660168320005\\
14.7106832971003	-0.227520067791415\\
15.1510454797550	-0.242441227098241\\
15.5914076624097	-0.258784289867716\\
16.0317698450644	-0.275931056891249\\
16.4721320277191	-0.293300510074889\\
16.9124942103738	-0.310362173917173\\
17.3528563930285	-0.326646684757160\\
17.7932185756832	-0.341753532540030\\
18.2335807583379	-0.355356026103084\\
18.6739429409926	-0.367203609275511\\
19.1143051236473	-0.377121720155885\\
19.5546673063020	-0.385009438978414\\
19.9950294889567	-0.390835210630198\\
20.4353916716114	-0.394630956164105\\
20.8757538542661	-0.396476344881239\\
21.3161160369208	-0.396502815018298\\
21.7564782195755	-0.394890009846882\\
22.1968404022302	-0.391846833417229\\
22.6372025848849	-0.387605235351901\\
23.0775647675396	-0.382408752810503\\
23.5179269501943	-0.376504997242476\\
23.9582891328490	-0.370137995180495\\
};

\addplot [color=black, dashed]
table[row sep=crcr]{%
0 -0.8\\
23.9582891328490 -0.8\\
};

\end{axis}
\end{tikzpicture}%
    \caption{Dynamic simulation from an example SUC solution.}
	\label{FigFrequencySimulation}
\end{figure}

In order to validate the frequency-security constraints obtained in Section \ref{SectionFrequency}, dynamic simulations of post-fault frequency evolution were run using MATLAB/Simulink. The system presented in Fig. \ref{FigBlockDiagram} was considered. The model used for generator dynamics is equivalent to that considered in \cite{OMalleyDeload,UCFaroe}, consisting of a droop control and a first-order model for governor dynamics, followed by a saturation block. The model used for BESS dynamics is based on the one considered in \cite{BESSdynamicsUK}, again using a first-order dynamics block.

An example solution from the proposed SUC model with binding nadir constraint is fed into the dynamic-model in Simulink. This SUC solution scheduled values of $H=132\textrm{GWs}$, $R_{\mathcal{S}}=0.22\textrm{GW}$, $R_{\mathcal{G}}=2.24\textrm{GW}$, $P_\textrm{L}=1.66\textrm{GW}$, and the demand level was $\textrm{P}_\textrm{D}=38.3\textrm{GW}$. Time-constants of $\tau_\textrm{g}=5\textrm{s}$ and $\tau_\textrm{b}=0.1\textrm{s}$ were considered for the simulation. As shown in Fig. \ref{FigFrequencySimulation}, the nadir obtained from the dynamic simulation is 0.72Hz, while the limit for nadir was set to $\Delta f_\textrm{max}=0.8\textrm{Hz}$ in the optimisation, demonstrating that the frequency nadir requirement is respected. The simulated RoCoF and q-s-s values are also within limits, $0.31\textrm{Hz/s}$ and $0.35\textrm{Hz}$ respectively. The conservativeness in nadir is driven by two factors: first, the analytical approximation in constraint (\ref{nadirEFRDamping}), which is further analysed in Section \ref{QuantitativeApproxSection}; and second, the approximation of droop control for FR by an increasing ramp of power injection in (\ref{EFRdefinition}) and (\ref{PFRdefinition}). This increasing-ramp assumption was demonstrated by \cite{OPFChavez} to conservatively approximate any generic droop control. Therefore, even if a more complex control of frequency dynamics than that shown in Fig. \ref{FigBlockDiagram} were to be considered, the proposed security constraints would still respect the requirements, as long as the regulation for EFR being delivered by $\textrm{T}_\textrm{s}$ and PFR being delivered by $\textrm{T}_\textrm{g}$ is complied with.

\subsection{Assessment of the Proposed Analytical Approximations} \label{QuantitativeApproxSection}

Three conservative approximations are proposed in this paper to linearize the nadir constraint (\ref{nadirEFRDamping}) while always guaranteeing frequency security: 1) a linear term approximating the effect of load damping on nadir, defined in (\ref{alphaEq}); 2) the overestimation by using planes for inner-approximating the squared term in the left-hand side of (\ref{nadirEFRDamping}), which is defined in (\ref{nadirPlanes}); and 3) the binary expansion of $R_\mathcal{G}$ as described in (\ref{binaryExpansion}). In this section, the conservativeness of these approximations is quantified, as well as their computational performance for different levels of accuracy of the approximations.

In order to assess the accuracy of the proposed linear damping term, inequality (\ref{nadirEFRDamping}) is used to obtain system conditions that exactly meet the nadir requirement. Then, for these conditions the actual value of nadir is computed by numerically solving the swing equation (\ref{SwingEq}). The same procedure is used for the case in which load damping is ignored, defined by constraint (\ref{nadirConstraint}). The difference between the actual frequency nadir and the nadir requirement provides an indication on the conservativeness of this approximation. For (\ref{nadirEFRDamping}) and (\ref{nadirConstraint}), 3500 samples of system conditions that exactly meet the inequalities were obtained, for system conditions in the ranges $H \in [50,400]\textrm{GW}\textrm{s}$, $R_{\mathcal{S}} \in [0,400]\textrm{MW}$, $R_{\mathcal{G}} \in [500,2500]\textrm{MW}$, $P_{\textrm{L}} \in [1400,1800]\textrm{MW}$ and $\textrm{P}_\textrm{D} \in [20,60]\textrm{GW}$, corresponding to GB's 2030 system.


\begin{table}[!t]
\renewcommand{\arraystretch}{1.35}
\caption{Effective Nadir Requirement}
\label{table_damping}
\centering
\begin{tabular}{|>{\centering\arraybackslash}m{1.4cm}|>{\centering\arraybackslash}m{3cm}|>{\centering\arraybackslash}m{2.85cm}|}
\cline{2-3}
    \multicolumn{1}{c|}{} & From constraint (\ref{nadirEFRDamping}), linear approx. of damping  &  From constraint (\ref{nadirConstraint}), neglecting damping \\ 
\hline
Mean nadir & 0.75Hz & 0.61Hz \\
\hline
Max. nadir & 0.79Hz & 0.72Hz \\
\hline
Min. nadir & 0.73Hz & 0.49Hz \\
\hline
\end{tabular}
\end{table}

\begin{table}[!t]
\renewcommand{\arraystretch}{1.4}
\caption{Assessment of Overestimating Planes' Approximation}
\label{table_planes}
\centering
\begin{tabular}{|>{\centering\arraybackslash}m{3.7cm}|>{\centering\arraybackslash}m{1.27cm}|>{\centering\arraybackslash}m{0.962cm}|>{\centering\arraybackslash}m{0.962cm}|}
\cline{2-4}
    \multicolumn{1}{c|}{} & 2 Planes (base case) & 4 Planes & 8 Planes \\ 
\hline
Mean overestimation in the entire range of $P_{\textrm{L}},R_{\mathcal{S}}$ & 0.72\% & 0.17\% & 0.04\% \\
\hline
Maximum overestimation at any point in the range of $P_{\textrm{L}},R_{\mathcal{S}}$ & 1.25\% & 0.35\% & 0.09\% \\
\hline
Increase in computational time 
& 0\% & 7\% & 60\% \\
\hline
Decrease in optimal cost of frequency services 
& 0\% & 0.4\% & 0.7\% \\
\hline
\end{tabular}
\end{table}

As shown in Table \ref{table_damping}, the proposed approximation on load damping leads to a 6\% conservativeness (0.05Hz) on average, while ignoring damping causes 25\% conservativeness (0.19Hz). Even in the worst case, the proposed approximation only increases the requirement by 9\% (0.07Hz), compared to 39\% (0.31Hz) when load damping is ignored. Table \ref{table_damping} demonstrates that the linear approximation is slightly conservative but very close to the actual nadir requirement of $\Delta f_\textrm{nadir}=0.8\textrm{Hz}$. Note that implementing constraint (\ref{nadirEFRDamping}) instead of (\ref{nadirConstraint}) does not affect the computational time of the optimisation problem. 

The conservativeness of the second approximation, i.e. the overestimation due to using planes for inner-approximating the squared term in (\ref{nadirEFRDamping}), is shown in Table \ref{table_planes} for different numbers of planes. These results were obtained again for $P_{\textrm{L}} \in [1400,1800]\textrm{MW}$ and $R_{\mathcal{S}} \in [0,400]\textrm{MW}$. There is a clear trend showing that using more planes reduces the conservativeness of the approximation but increases the computational time of the SUC. By increasing the number of planes from 2 to 4 and then to 8, the cost of frequency services (calculated after running the SUC as the total system operating cost net the energy cost) reduces by 0.4\% and 0.7\%, while the computational time is increased by 7\% and 60\%, respectively.

Regarding the binary expansion of $R_\mathcal{G}$, the precision is determined by the number of bits used. For the purpose of this paper, a set of 12 bits, $\mathcal{L}=\{0,1,...\,,11\}$, is used as the base case, since the highest representable value of 4095MW is an acceptable upper bound for $R_\mathcal{G}$ in GB's 2030 system. From this base case, some of the Least Significant Bits (LSBs) are removed, in order to study the improvement in computational burden due to renouncing to some precision in the binary expansion. The results in Table \ref{table_bits} clearly show that a lower number of binary decision variables significantly reduces computational time for the optimisation with the drawback of increased conservativeness. 

Since Stochastic Programming is a computationally demanding problem, a balance between accuracy in the optimal objective and computational efficiency is desired for the SUC-based simulations presented in coming sections. Therefore, 2 planes and 7 bits (5 LSBs removed) for the expansion of $R_{\mathcal{G}}$ are considered from now on. For future applications of these frequency-security constraints, different numbers of planes and bits can be chosen, depending on the time-requirements for solving the optimisation in each application.

\begin{table}[!t]
\renewcommand{\arraystretch}{1.4}
\caption{Assessment of Approximation Due to Binary Expansion of $R_\mathcal{G}$}
\label{table_bits}
\centering
\begin{tabular}{|>{\centering\arraybackslash}m{3.0cm}|>{\centering\arraybackslash}m{1.27cm}|>{\centering\arraybackslash}m{1.17cm}|>{\centering\arraybackslash}m{1.17cm}|}
\cline{2-4}
    \multicolumn{1}{c|}{} & 12 bits (base case) & 5 LSBs removed & 10 LSBs removed \\ 
\hline
Precision achieved & 1MW & 32MW & 1024MW \\
\hline
Decrease in computational time 
& 0\% & 70\% & 97\% \\
\hline
Increase in optimal cost of frequency services 
& 0\% & 1\% & 49\% \\
\hline
\end{tabular}
\end{table}

\subsection{Value of Defining and Optimising EFR as a distinct service} \label{EFRvalue}

This section presents the benefits from recognising the faster dynamics of EFR and  optimising its provision in the scheduling process. Three different operating strategies for GB's power system are considered: 1) ``Just PFR" (considered as the base-case strategy): EFR is not defined as a distinct service, as was the traditional practice. Therefore, the FR provided by BESS is considered as PFR; 2) ``Fixed EFR": EFR is defined as service, but its provision is fixed throughout the year, reflecting the current practice in GB; 3) ``Optimised EFR": EFR is defined and co-optimised along with PFR and inertia, using the model proposed in this paper. For all these strategies, the largest power infeed is kept constant to $P_{\textrm{L}}=\textrm{P}^{\textrm{max}}_{\textrm{L}}$. 

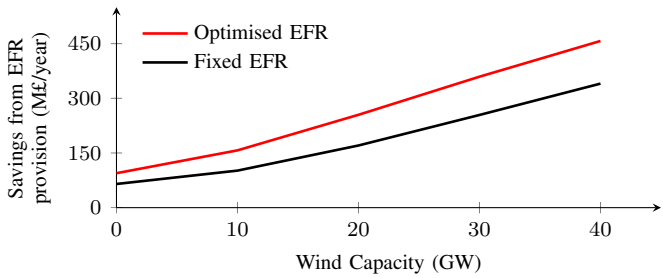
\begin{figure}
\hspace*{-0.2cm}
%
%
\definecolor{mycolor1}{rgb}{1.00000,0.26275,0.26275}%
\begin{tikzpicture}

\begin{axis}[%
axis lines = left,
width=2.85in,
height=1.05in,
at={(0in,0in)},
scale only axis,
y label style={at={(axis description cs:-0.16,-0.03)},anchor=west},
xmin=0,
xmax=45,
xlabel style={font=\footnotesize},
xlabel={Wind Capacity (GW)},
ymin=0,
ymax=550,
ytick={0,150,300,450},
ylabel style={font=\footnotesize},
ylabel style={align=center},
ylabel={Savings from EFR\\provision (M\pounds/year)},
axis background/.style={fill=white},
legend style={at={(0.03,0.97)}, anchor=north west, legend cell align=left, align=left, draw=white!15!black},legend style={font=\footnotesize},
legend style={fill opacity=0,text opacity=1,draw=none}, 
xticklabel style={font=\footnotesize},
yticklabel style={font=\footnotesize}
]




\addplot [color=red, line width=1.0pt]
  table[row sep=crcr]{%
0	94.51164\\
10 157.1894\\
20	254.7846\\
30 359.3352\\
40	457.32456\\
};
\addlegendentry{Optimised EFR}

\addplot [color=black, line width=1.0pt]
  table[row sep=crcr]{%
0	64.93788\\
10 101.6248\\
20	170.54844\\
30 254.3116\\
40	340.34352\\
};
\addlegendentry{Fixed EFR}


\end{axis}
\end{tikzpicture}%
    \caption{Annualised savings due to EFR provision from a 200MW-rated BESS, for two different operating strategies. Both strategies' savings are referred to the case in which the BESS provides just PFR.}
	\label{FigEFRWindSavings}
\end{figure}

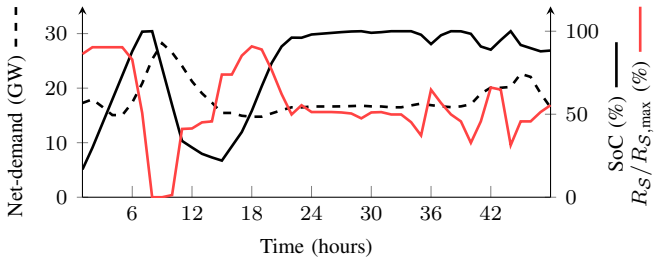
\begin{figure}
\hspace*{-0.2cm}
\begin{tikzpicture}
\pgfplotsset{set layers}
\definecolor{mycolor1}{rgb}{1.00000,0.26275,0.26275}%
\begin{axis}[
axis y line=left,
axis x line*=left,
width=2.45in,
height=1in,
at={(0in,0in)},
scale only axis,
xmin=1,
xmax=48,
xtick={0,6,12,18,24,30,36,42},
xlabel style={font=\footnotesize},
xlabel={Time (hours)},
ymin=0,
ymax=35,
ylabel style={font=\footnotesize},
ylabel style={align=center},
y label style={at={(axis description cs:-0.14,-0.15)},anchor=west},
ylabel={Net-demand (GW) \ref{pgfplots:plot1_2}},
xticklabel style={font=\footnotesize},
yticklabel style={font=\footnotesize}
]
\addplot [color=black, dashed, line width=1.0pt]
  table[row sep=crcr]{%
1	17.27648\\
2	17.94454\\
3	16.57624\\
4	15.06554\\
5	15.05188\\
6	17.36693\\
7	20.87004\\
8	25.29356\\
9	28.22265\\
10	26.52667\\
11	24.0721\\
12	21.24771\\
13	19.05512\\
14	17.34999\\
15	15.454232671\\
16	15.454223\\
17	14.93999871\\
18	14.7598865\\
19	14.7688796\\
20	15.371995\\
21	15.95963\\
22	16.487927\\
23	16.465664\\
24	16.6383\\
25	16.63834\\
26	16.638365\\
27	16.65883\\
28	16.71028\\
29	16.78141\\
30	16.658758\\
31	16.640377\\
32	16.487358\\
33	16.487358\\
34	16.83369\\
35	17.17106469\\
36	16.87832\\
37	16.66184\\
38	16.5155\\
39	16.65711391\\
40	17.04119\\
41	18.98637\\
42	20.09952\\
43	20.06999\\
44	20.22221\\
45	22.62542\\
46	22.06364\\
47	19.053\\
48	16.4603923\\
};
\label{pgfplots:plot1_2}
\end{axis}
\begin{axis}[
axis y line=right,
axis x line=none,
width=2.45in,
height=1in,
scale only axis,
xmin=1,
xmax=48,
ymin=0,
ymax=115,
ylabel style={font=\footnotesize},
ylabel style={align=center},
ylabel={SoC (\%) \ref{pgfplots:plot2_2} \\ $R_{\mathcal{S}}/R_{\mathcal{S},\textrm{max}}$ (\%) \ref{pgfplots:plot3_2}},
y label style={at={(axis description cs:0.985,-0.1)},anchor=west},
yticklabel style={font=\footnotesize}
]

\addplot [color=black, line width=1.0pt]
table[row sep=crcr]{%
1	16.71774\\
2	29.78697\\
3	44.31289\\
4	58.83882\\
5	73.36474\\
6	87.89066\\
7	99.74041\\
8	100\\
9	77.77778\\
10	55.55556\\
11	33.9148\\
12	29.97539\\
13	26.16053\\
14	23.98001\\
15	22.06149\\
16	30.66729\\
17	39.27309\\
18	51.97692\\
19	66.82365\\
20	80.86116\\
21	90.82952\\
22	96.13691\\
23	96.07353\\
24	98.02974\\
25	98.50646\\
26	98.98216\\
27	99.45775\\
28	99.84275\\
29	100\\
30	99.04657\\
31	99.43189\\
32	100\\
33	100\\
34	100\\
35	97.87113\\
36	92.25176\\
37	97.55705\\
38	100\\
39	99.94829\\
40	98.44595\\
41	90.79694\\
42	88.87842\\
43	94.66555\\
44	100\\
45	91.72167\\
46	89.80315\\
47	87.88459\\
48	88.3613\\
};
\label{pgfplots:plot2_2}

\addplot [color=mycolor1, line width=1.0pt]
  table[row sep=crcr]{%
1 86.3034250000000\\
2 90.3497750000000\\
3 90.3497750000000\\
4 90.3497750000000\\
5 90.3497750000000\\
6 82.915950000000\\
7 50.7211000000000\\
8 0\\
9 0\\
10 01.30579775000000\\
11 41.1363250000000\\
12 41.4165750000000\\
13 45.0938250000000\\
14 45.6808250000000\\
15 73.9025000000000\\
16 73.9025000000000\\
17 85.2884250000000\\
18 90.8665750000000\\
19 88.9906000000000\\
20 76.6787750000000\\
21 62.0918000000000\\
22 49.8548750000000\\
23 55.4339250000000\\
24 51.3217250000000\\
25 51.3214000000000\\
26 51.3210750000000\\
27 51.0694250000000\\
28 50.4368250000000\\
29 47.5152000000000\\
30 51.0703500000000\\
31 51.2817250000000\\
32 49.8548750000000\\
33 49.8548750000000\\
34 45.2100500000000\\
35 37.2074250000000\\
36 64.7369000000000\\
37 56.7859750000000\\
38 49.8811500000000\\
39 45.6808250000000\\
40 32.7872250000000\\
41 45.6808250000000\\
42 66.0753500000000\\
43 64.8179250000000\\
44 31.3737500000000\\
45 45.6808250000000\\
46 45.6808250000000\\
47 51.3217250000000\\
48 55.4345750000000\\
};
\label{pgfplots:plot3_2}

\end{axis}
\end{tikzpicture}
    \caption{Two-day example of the operation of GB's 2030 system.}
	\label{FigOperation}
\end{figure}

Fig. \ref{FigEFRWindSavings} displays the operational cost savings from strategies 2 and 3 referred to the base-case strategy 1, therefore showing the value of defining EFR as a distinct service and optimising its provision, respectively. BESS of 200MW rating and 5h tank size are considered, for five cases of wind capacity corresponding to 0, 10, 20, 30 and 40GW. As expected, the value of EFR increases with increasing wind penetration, as the declining system inertia makes EFR more valuable. Furthermore, the savings from applying the proposed model to optimise EFR provision are higher than 33\% for any level of wind penetration, when compared to just providing a fixed amount. The results demonstrate the need to not only recognise the faster dynamics of EFR but also optimise its provision. The higher benefits from optimising EFR provision come from two aspects: by optimising EFR, a 200MW-BESS can provide up to 400MW of EFR at times; in addition the BESS can provide other services such as reserve. These two aspects are better understood by looking at the detailed operation of the BESS, shown in Fig. \ref{FigOperation}.

Fig. \ref{FigOperation} presents a two-day period operation of BESS using strategy 3 : it displays the net-demand, the State of Charge (SoC) of the BESS, and the EFR scheduled as a percentage of the maximum EFR that could be provided, $R_{\mathcal{S},\textrm{max}}$. Note that $R_{\mathcal{S},\textrm{max}}$ takes the value of double the rating of the BESS, therefore $R_{\mathcal{S},\textrm{max}}=400\textrm{MW}$ in this case; indeed an EFR of up to 400MW can be provided by BESS by swiftly shifting from fully charging to fully discharging, which would effectively increase power injection by twice the volume of the BESS. Therefore, an optimised provision of EFR allows BESS to provide a higher amount of EFR in periods of low net-demand (and therefore low inertia) as hours 4 and 18, which is achieved by scheduling the BESS to fully charge. On the other hand, when the system inertia is high (e.g. around hour 8), the BESS provides zero EFR while discharging. 

Furthermore, the ability to shift between charging/discharging modes, when compared to the ``Fixed EFR'' strategy in which the BESS would be forced to stay idle in order to provide 200MW of EFR, also allows to take advantage of a synergy with energy costs. For low net-demand periods, EFR is more valuable while the cost of energy is generally lower, so charging the BESS in these periods is positive both from EFR and energy perspectives. In addition, the BESS can provide reserve when needed.

\begin{figure}
%
%
\definecolor{mycolor1}{rgb}{1.00000,0.26275,0.26275}%
\begin{tikzpicture}

\begin{axis}[%
axis lines = left,
width=2.65in,
height=1in,
at={(0in,0in)},
scale only axis,
y label style={at={(axis description cs:-0.16,0)},anchor=west},
xmin=0,
xmax=950,
xlabel style={font=\footnotesize},
xlabel={BESS volume (MW)},
xtick={0,100,200,300,400,500,600,700,800,900},
ymin=0,
ymax=48,
ylabel style={font=\footnotesize},
ylabel style={align=center},
ylabel={Savings increase\\(\%/year)},
axis background/.style={fill=white},
legend style={at={(0.03,0.97)}, anchor=north west, legend cell align=left, align=left, draw=white!15!black},legend style={font=\footnotesize},
legend style={fill opacity=0,text opacity=1,draw=none}, 
xticklabel style={font=\footnotesize},
yticklabel=\pgfmathprintnumber{\tick}\%,
yticklabel style={font=\footnotesize}
]






\addplot [color=black, line width=1.0pt]
  table[row sep=crcr]{%
50 40.1946\\
100 34.3365\\
200	33.792238912\\
400	26.114181746\\
600	24.703222611\\
750 22.7932\\
900	08.550521965\\
};

\end{axis}
\end{tikzpicture}%
    \caption{Increase in savings from strategy ``Optimised EFR" with respect to strategy ``Fixed EFR", for a 40GW-wind scenario.}
	\label{FigEFRBESSSavings}
\end{figure}
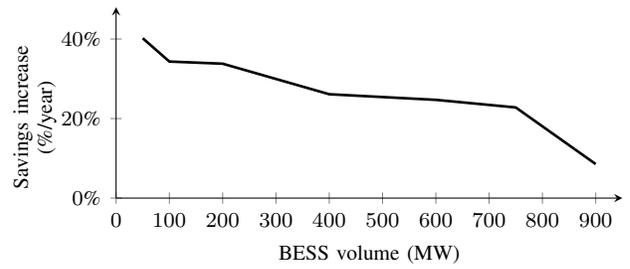

Finally, Fig. \ref{FigEFRBESSSavings} presents the extra savings from ``Optimised EFR" provision over the ``Fixed EFR" case. The volume of BESS is increased up to 900MW. These results show that the benefit from optimising EFR provision is very significant when a small volume of BESS is available, but reduces as the volume of BESS increases. For cases with more than 800MW BESS available, the provision of EFR is sufficient even for periods of low net-demand, leading to very limited benefit of optimising its provision. 

\subsection{Value of Dynamically-Reduced Largest Power Infeed}

\begin{figure}[!t]
\centering
\includegraphics[width=3.4in]{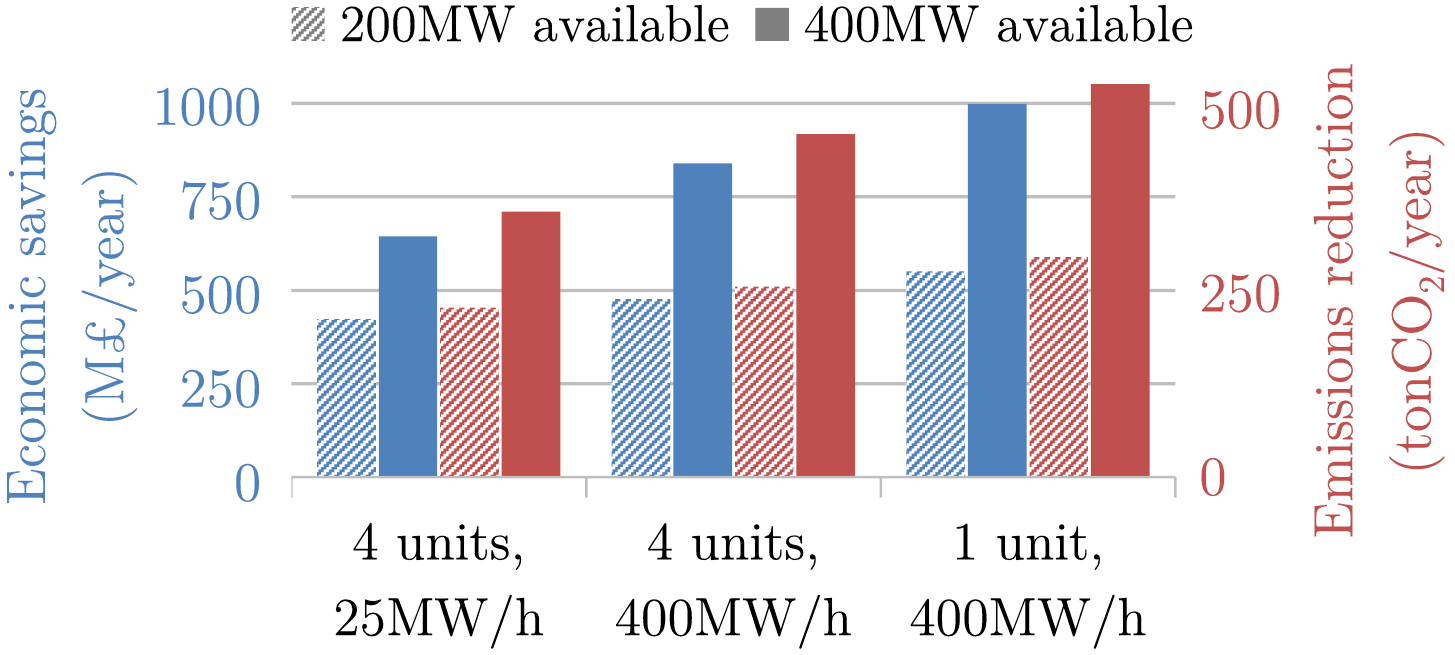}
\caption{Economic and emissions savings due to nuclear deloading, for deloading availabilities of 200MW and 400MW in a 40GW-wind case.}
\label{BarChart}
\end{figure}

A lower largest-power-infeed, represented by a decreasing value of decision variable $P_\textrm{L}$ in this paper, effectively reduces the need for frequency services, as can be observed in the frequency-security constraints presented in Section \ref{SectionFrequency}. Although nuclear units provide low-cost, zero-emissions energy, the large capacity of each single unit becomes the key driver of the frequency-response requirement in GB. Under certain system's conditions, it may be cost-effective to reduce their power output in order to reduce the need for frequency services. Note that deloading nuclear units would be used as a preventive measure only when a low net-demand period is expected, since ramping capabilities of these types of units are typically limited, and therefore must be deloaded some time in advance. 

Fig. \ref{BarChart} shows the operating cost savings from a dynamically-reduced largest power infeed, with respect to the ``Just PFR" strategy discussed in Section \ref{EFRvalue}. A 40GW wind scenario is considered. Two different deloading capabilities are studied: 200MW and 400MW. The impact of the nuclear fleet's characteristics is analysed by considering 4 large units rated at 1.8GW and just 1 large unit. 
In addition, two extreme ramping rates of the units are considered, 25MW/h and 400MW/h. The results suggest that by allowing the deloading of nuclear units, considerable cost saving can be achieved. However, a higher number of large units and lower ramp rates reduce the savings, due to the fact that all large nuclear plants need to be deloaded in order to reduce the largest infeed, and the fact that it is more challenging for slower plants to lower their power outputs when required.

Furthermore, deloading nuclear units can reduce carbon emissions. Although this might seem contradictory, lowering the power production of large carbon-free nuclear plants reduces the need for frequency services and consequently less part-loaded conventional generators need to be online, allowing for more RES to be accommodated. 

Note that the load factor of nuclear units would be reduced if the deloading strategy was adopted. Subsidies may need to be introduced to compensate for the reduced power production, which would reduce the savings from this strategy.

\subsection{Full Co-Optimisation of Frequency Services} \label{SectionFullOpt}
\begin{figure}[!t] 
\centering
\includegraphics[width=3.3in]{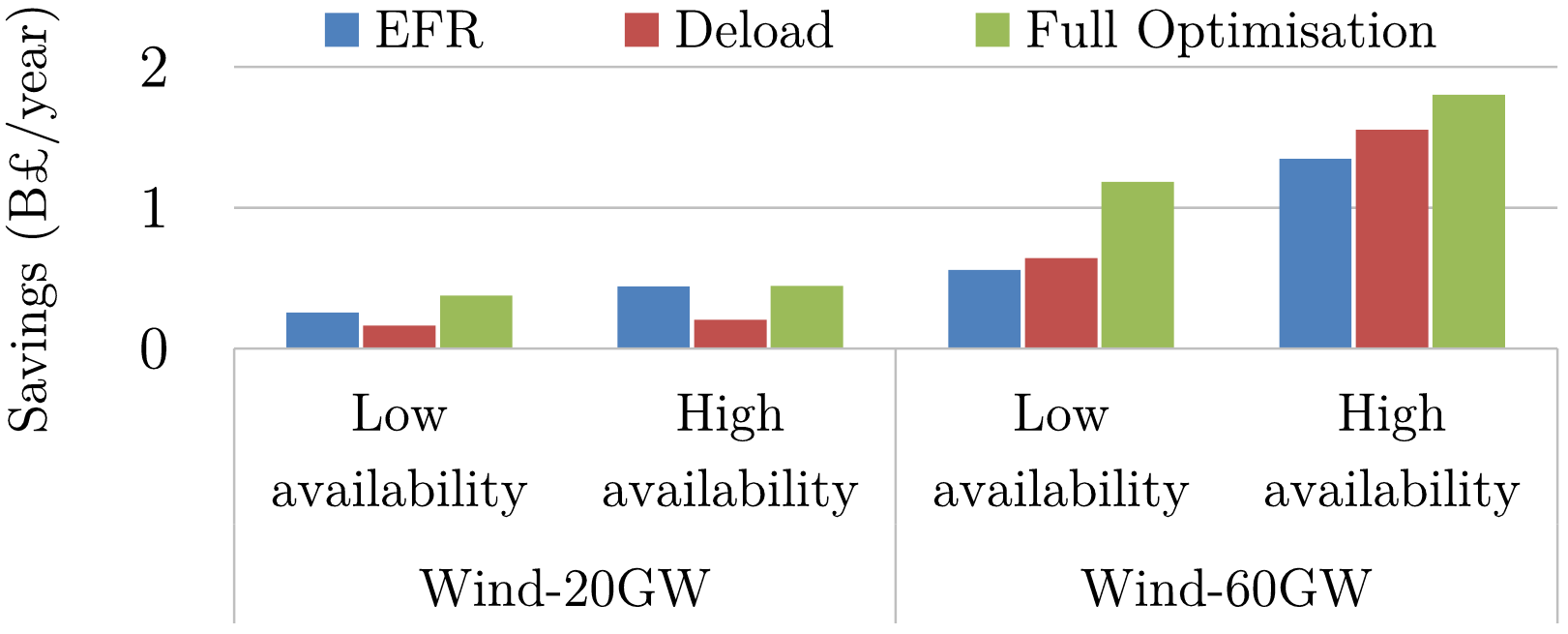}
\caption{Annualised savings from optimising frequency services.  The \textit{EFR} strategy presents the savings from an optimised EFR provision, with no option to deload nuclear plants; the \textit{Deload} strategy presents the savings from optimising the deloading option of large nuclear units, with no EFR availability in the system; finally, \textit{Full Optimisation} considers the co-optimisation of both these services.}
\label{FinalPlot}
\end{figure}

The benefits of simultaneously optimising all four frequency services in the SUC are presented here. Two cases are considered: a ``Low availability" case, corresponding to a 200MW-rated BESS and a maximum deloading of 200MW; and a ``High availability" case, corresponding to a 600MW-rated BESS and a maximum deloading of 600MW. The largest power infeed is driven by 4 nuclear units rated at $\textrm{P}_\textrm{L}^\textrm{max}=1.8\textrm{GW}$, with 100MW/h ramp rates. The results for both cases are presented in Fig. \ref{FinalPlot}, showing the cost savings from each case referred to the ``Just PFR" strategy discussed in Section \ref{EFRvalue}. The savings from this ``Full optimisation'' are compared to the savings from just EFR optimisation and just deloading optimisation.

EFR shows to be more beneficial than deloading of nuclear units for the low-wind scenario, as EFR provision is free of cost in this framework but deloading nuclear units increases the cost of energy in low-wind scenarios, since the energy not produced by nuclear must be provided by more expensive thermal plants. However, for high wind penetration deloading becomes more valuable. In this scenario, the energy not provided from nuclear plants due to deloading is now provided by wind, free of cost. In addition, deloading nuclear plants is more valuable for securing post-fault frequency, as its delivery time is virtually zero while EFR takes $\textrm{T}_\textrm{s}$ to be delivered.

Furthermore, the savings from the ``Full Optimisation" strategy are significantly lower than the sum of savings from just EFR and just deloading optimisations, in the ``High availability" cases. This result shows clear competition between these two services if there is a significant amount of both services available. 


\subsection{Impact of Damping and SI in the Value of Frequency Services} \label{sectionSI}

\begin{figure}[!t] 
\centering
\includegraphics[width=3.3in]{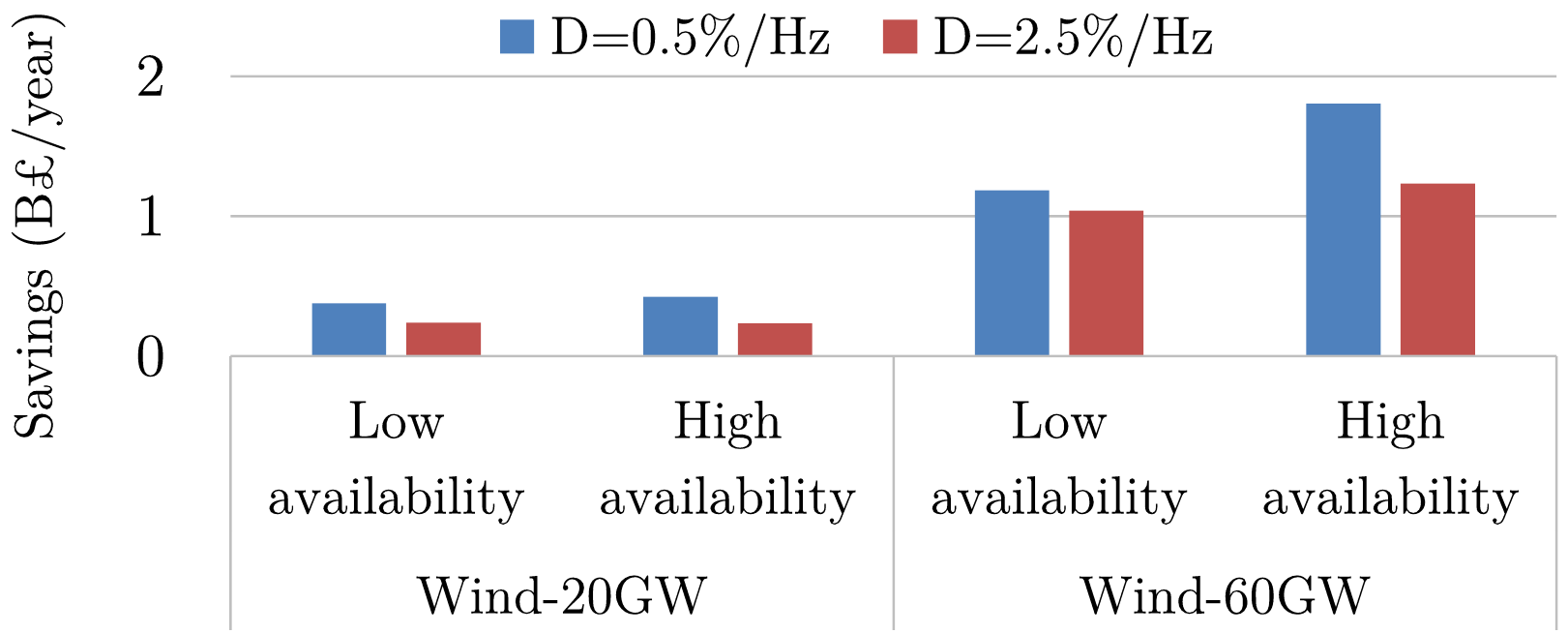}
\caption{Impact of load damping on the savings from frequency services, for a \textit{Full Optimisation} strategy.}
\label{FinalPlot2}
\end{figure}

This section analyses the reduction of the value of EFR and nuclear deloading in the presence of higher load damping and synthetic inertia provision from wind turbines. Fig. \ref{FinalPlot2} presents the impact of load damping on the benefits from introducing these two new frequency services. The savings from EFR and deloading are not significantly affected by higher damping when a limited amount of these services is available, but the savings become much more sensitive to damping if the availability of the services increases. This result suggests that increasing system damping would not affect the competitiveness of EFR and deloading but would instead reduce the required volume of services.

\begin{figure}[!t] 
\centering
\includegraphics[width=3.3in]{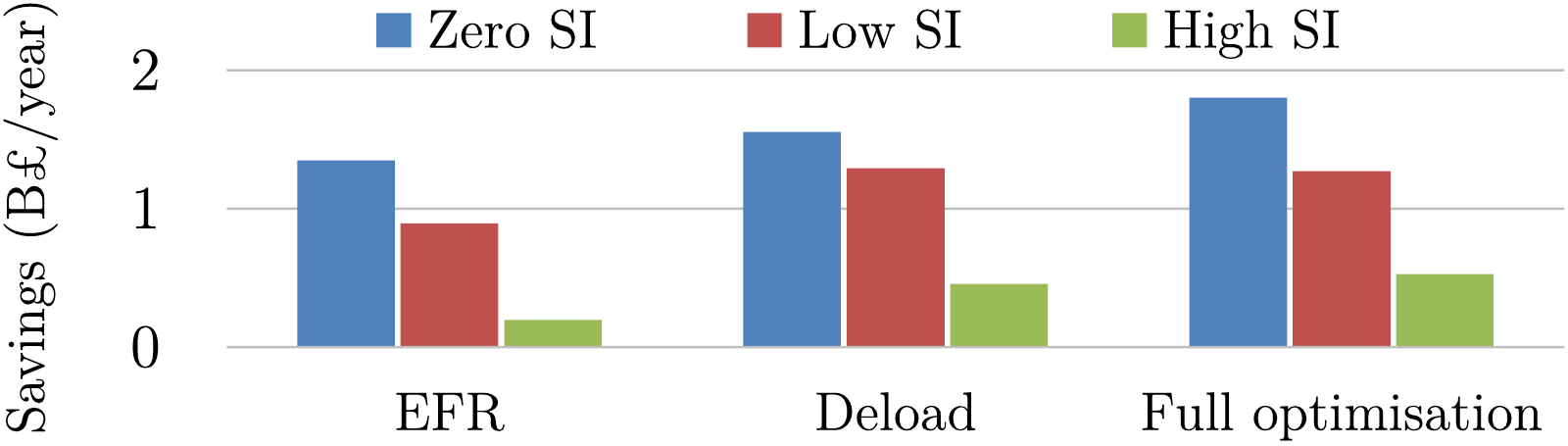}
\caption{Impact of different levels of synthetic inertia (SI) on the value of other frequency services, for a 60-GW-wind-capacity case and ``High availability'' of EFR and nuclear deloading.}
\label{FinalPlot3}
\end{figure}

Finally, Fig. \ref{FinalPlot3} displays the impact of SI provision from wind turbines on the value of other frequency services. The model developed in \cite{FeiAssessment} for considering SI in a scheduling algorithm was integrated in the frequency-secured SUC proposed in this paper. One of the main difficulties in modelling SI in the system level comes from the uncertainty in the number of wind turbines online for a given level of wind power generation. This paper has assumed the average number of online turbines as in \cite{SyntheticOMalley}, while considering two different inertia time-constants for individual wind turbines: in ``Low SI'' $\textrm{H}_\textrm{W}=1.5\textrm{s}$, while in ``High SI'' $\textrm{H}_\textrm{W}=5\textrm{s}$. The results in Fig. \ref{FinalPlot3}, for a high wind-penetration case of 60GW and the ``High availability'' case discussed in Section \ref{SectionFullOpt}, suggest that the potential development of SI capability from RES needs to be taken into consideration when designing alternative frequency services. In the case with limited SI capability from wind turbines, there are still clear benefits from introducing the EFR and nuclear-deloading services, while in the case that wind turbines have a similar inertia time-constant to conventional plants, the benefits and need for alternative frequency services are very limited.

\subsection{Comparison of Results in Stochastic and Deterministic UC} \label{SectionDetVsStoch}

\begin{table}[!t]
\renewcommand{\arraystretch}{1.35}
\caption{Comparison of Results, from Deterministic UC\hspace{\textwidth}over Stochastic UC}
\label{table_DetStoch}
\centering
\begin{tabular}{|>{\centering\arraybackslash}m{3cm}|>{\centering\arraybackslash}m{2cm}|>{\centering\arraybackslash}m{2cm}|}
\cline{2-3}
    \multicolumn{1}{c|}{} & Wind Capacity 20GW  &  Wind Capacity 60GW \\ 
\hline
Total Operating Cost & 0.6\% & 1\% \\
\hline
Wind Curtailment & 0.9\% & 1.4\% \\
\hline
Savings from EFR and Deload & -12\% & -15\% \\
\hline
Average $H$ & 5\% & 7\% \\
\hline
Average $R_\mathcal{S}$ & -4\% & -5\% \\
\hline
Average $P_\textrm{L}$ & 1\% & 2\% \\
\hline
\end{tabular}
\end{table}

Here we compare the proposed frequency-constrained framework in both the SUC and a Deterministic UC. The ``Full Optimisation" strategy presented before was used in a ``Low Availability" case for frequency services (200MW-rated BESS and a maximum nuclear deloading of 200MW), while two wind-capacity cases were considered, of 20GW and 60GW. The Deterministic UC was run using a quantile of 0.98 for net-demand, as in \cite{AlexEfficient}. The comparison of results from both UC approaches are presented in Table \ref{table_DetStoch}. This table shows the difference, in percentage, of several system magnitudes from the solution of the Deterministic UC with respect to the solution of the SUC. 

Table \ref{table_DetStoch} demonstrates that if the Deterministic UC is applied, the system operating cost and wind curtailment are both higher than for SUC. At the same time, the savings from the new frequency services, namely EFR and part-loading nuclear, are lower in Deterministic UC when compared to SUC. This is due to the fact that Deterministic UC tends to schedule a higher number of slow-start thermal generators (CCGTs) as shown in \cite{AlexEfficient}, implying a higher level of inertia than in SUC. The higher inertia makes EFR and Deloading less valuable to comply with the frequency constraints. Nevertheless, the savings obtained from co-optimising these services are still significant under a Deterministic UC: \pounds320m/year for the 20GW-wind case and \pounds980m/year for the 60GW-wind case.

\section{Conclusion and Future Work} \label{SectionConclusion}
This paper proposes an SUC framework for co-optimising energy production under uncertainty, along with the scheduling of diverse frequency services, namely synchronised and synthetic inertia, Primary Frequency Response (PFR), EFR and a dynamically-reduced largest power infeed.
A set of linear constraints are developed, which guarantee the fulfillment of post-fault frequency requirements. Their accuracy and computational efficiency for implementation in the very computationally-demanding SUC problem has been investigated. Several case studies clearly demonstrate the benefits of introducing new frequency response services. The results also highlight the importance of co-optimisation of alternative services to accurately capture their value.

Three main areas for enhancing the proposed model can be identified. Firstly, the optimisation of several FR products with distinct delivery times, other than just EFR and PFR, would allow to recognise the different dynamics of providers, leading to a more efficient ancillary-services market. Secondly, the present paper considers the uniform frequency model, while diverse frequency evolutions across the network's buses have been recognised. Further work needs to be carried out to investigate the impact of locational frequencies on the value of these frequency services. Thirdly, the impact of frequency-security constraints on pricing schemes would be of interest for the design of an efficient ancillary-services market for inertia and FR, so the work in \cite{ElaI,PricingElaZhang} could be enhanced to incorporate the new frequency services introduced in the present paper.

%
%
%
%



\ifCLASSOPTIONcaptionsoff
  \newpage
\fi



%



\bibliographystyle{IEEEtran} 
\bibliography{Luis_PhD}




\newpage

\begin{IEEEbiographynophoto}{Luis Badesa}
(S'14) received the B.S. in Industrial Engineering degree 
from the University of Zaragoza, Spain, in 2014, and the 
M.S. in Electrical Engineering degree from the University
of Maine, United States, in 2016. Currently he is pursuing
a Ph.D. in Electrical Engineering at Imperial 
College London, U.K. His research interests lie in modelling
and optimisation of low-carbon power grids' operation.
\end{IEEEbiographynophoto}

 


\vspace{-13.2cm}

\begin{IEEEbiographynophoto}{Fei Teng}
(M'15) received the BEng in Electrical Engineering
from Beihang University, China, in 2009, and the 
MSc and PhD degrees in Electrical Engineering from 
Imperial College London, U.K., in 2010 and 2015.
Currently, he is a Lecturer in the Department of Electrical 
and Electronic Engineering, Imperial College London, U.K.
His research focuses on advanced modelling, optimisation 
and analysis of low-carbon power systems' operation and planning.
\end{IEEEbiographynophoto}


\vspace{-13.2cm}

\begin{IEEEbiographynophoto}{Goran Strbac}
(M'95) is a Professor of Electrical Energy
Systems at Imperial College London, U.K. His current
research is focused on optimisation of energy systems'
operation and investment, energy infrastructure reliability and 
future energy markets. 
\end{IEEEbiographynophoto}




\end{document}